\documentclass[paper=letter, DIV=12, numbers=noenddot]{scrartcl} 
\usepackage{cmap}		
\usepackage[utf8]{inputenc}	
\usepackage[T1]{fontenc}	
\usepackage[USenglish]{babel}
\usepackage[]{microtype}

\usepackage{amsmath}
\usepackage{amssymb}
\usepackage{mathtools}		
\usepackage{xfrac}

\usepackage{aliascnt}		
\usepackage{amsthm}

\usepackage{tikz}		
\usetikzlibrary{matrix,arrows,patterns,intersections,calc,decorations.pathmorphing}

\usepackage[pdftitle={Expected covering radius of a translation surface},pdfauthor={Howard Masur, Kasra Rafi, Anja Randecker},pdfborder={0 0 0}]{hyperref}   
\usepackage[figure, table]{hypcap}	

\newtheoremstyle{break}
 {} 
 {} 
 {\itshape} 
 {} 
 {\bfseries} 
 {} 
 {\newline} 
 {\thmname{#1}\thmnumber{ #2}\thmnote{ (#3)}} 
 
\newtheoremstyle{breakdef}
 {} 
 {} 
 {} 
 {} 
 {\bfseries} 
 {} 
 {\newline} 
 {\thmname{#1}\thmnumber{ #2}\thmnote{ (#3)}} 
 
\newtheoremstyle{remark}
 {} 
 {} 
 {} 
 {} 
 {\itshape} 
 {.} 
 {0.5em} 
 {\thmname{#1}\thmnumber{ #2}\thmnote{ {\normalfont (}#3{\normalfont )}}} 

\theoremstyle{breakdef}
\newtheorem{definition}{Definition}[section]

\theoremstyle{remark}

\newaliascnt{rem}{definition}  
\newtheorem{rem}[rem]{Remark}
\aliascntresetthe{rem}

\newaliascnt{exa}{definition}  

\aliascntresetthe{exa} 

\newaliascnt{lem}{definition}  
\newtheorem{lem}[lem]{Lemma}
\aliascntresetthe{lem}

\theoremstyle{break}

\newaliascnt{thm}{definition}  
\newtheorem{thm}[thm]{Theorem}
\aliascntresetthe{thm}

\newaliascnt{prop}{definition}  

\aliascntresetthe{prop}

\newaliascnt{cor}{definition}  
\newtheorem{cor}[cor]{Corollary}
\aliascntresetthe{cor}  

  \newtheorem{introthm}{Theorem}
  \newtheorem{introprop}[introthm]{Proposition}

\makeatletter
\newtheorem*{rep@theorem}{\rep@title}
\newcommand{\newreptheorem}[2]{%
\newenvironment{rep#1}[1]{%
 \def\rep@title{#2 \ref{##1}}%
 \begin{rep@theorem}}%
 {\end{rep@theorem}}}
\makeatother

\newreptheorem{cor}{Corollary}

\numberwithin{figure}{section}			

\newcommand{\cB}{{\mathcal{B}}}
\newcommand{\cC}{{\mathcal{C}}}
\newcommand{\cE}{{\mathcal{E}}}
\newcommand{\cH}{{\mathcal{H}}}
\newcommand{\cO}{{\mathcal{O}}}
\newcommand{\cM}{{\mathcal{M}}}

\newcommand{\Hsd}{{\mathcal{H}_{\mathrm{small\text{-}diam}}}}
\newcommand{\Htc}{{\mathcal{H}_{\mathrm{thin\text{-}cyl}}}}
\newcommand{\Hpc}{{\mathcal{H}_{\mathrm{poor\text{-}cyl}}}}
\newcommand{\Hrc}{{\mathcal{H}_{\mathrm{rich\text{-}cyl}}}}
\newcommand{\Hrd}{{\mathcal{H}_{\mathrm{rich\text{-}disk}}}}
\newcommand{\Hrdc}{{\mathcal{H}_{\leq 1, \mathrm{rich\text{-}disk}}}}
\newcommand{\Hrdu}{{\mathcal{H}_{\geq \frac{1}{2}, \mathrm{rich\text{-}disk}}}}
\newcommand{\Ht}{{\mathcal{H}_{\mathrm{thin}}}}

\DeclareMathOperator{\SL}{SL}

\DeclareMathOperator{\area}{area}
\DeclareMathOperator{\diam}{diam}
\DeclareMathOperator{\crad}{c\text{-}rad}
\DeclareMathOperator{\height}{height}
\DeclareMathOperator{\Jac}{Jac}
\DeclareMathOperator{\im}{im}
\DeclareMathOperator{\Vol}{Vol}

  \newcommand{\param}{{\mathchoice{\mkern1mu\mbox{\raise2.2pt\hbox{$
  \centerdot$}}
  \mkern1mu}{\mkern1mu\mbox{\raise2.2pt\hbox{$\centerdot$}}\mkern1mu}{
  \mkern1.5mu\centerdot\mkern1.5mu}{\mkern1.5mu\centerdot\mkern1.5mu}}}

\usepackage{enumerate}

\newcommand\blfootnote[1]{
  \begingroup
  \renewcommand\thefootnote{}\footnote{#1}
  \addtocounter{footnote}{-1}
  \endgroup
}

\author{
Howard Masur
\thanks{Department of Mathematics, University of Chicago, Chicago, Il.,
{\texttt{masur@math.uchicago.edu}}}
\and Kasra Rafi
\thanks{Department of Mathematics, University of Toronto, Toronto, ON,
{\texttt{rafi@math.toronto.edu}}}
\and Anja Randecker
\thanks{Department of Mathematics, University of Toronto, Toronto, ON,
{\texttt{anja@math.toronto.edu}}}
}

\title{Expected covering radius of \\ a translation surface}
\date{\today}

\begin{document}
\renewcommand{\sectionautorefname}{Section}
\renewcommand{\subsectionautorefname}{Subsection}

\maketitle

\begin{abstract}
A translation structure equips a Riemann surface with a singular flat metric. 
Not much is known about the shape of a random translation surface. 
We compute an upper bound on the expected value of the covering radius of a 
translation surface in any stratum $\cH_1(\kappa)$.
The covering radius of a translation surface is the largest radius of an immersed disk. 
In the case of the stratum $\cH_1(2g-2)$ of translation surfaces of genus~$g$ with 
one singularity, the covering radius is comparable to the diameter.  
We show that the expected covering radius of a surface is bounded above by 
a uniform multiple of~$\sqrt{ \frac {\log g}{g}}$, independent of the stratum. 
This is smaller than what one would expect 
by analogy from the result of Mirzakhani about the expected diameter of a hyperbolic 
metric on a Riemann surface. 
To prove our result, we need an estimate for the volume of the thin part of $\cH_1(\kappa)$ which is given in the appendix. 
\end{abstract} 

\blfootnote{\textup{2010} \textit{Mathematics Subject Classification}: \textup{32G15, 30F60, 57M50, 37P45}}

\section{Introduction}

Translation surfaces have been studied in depth for many years. However, 
there is no clear picture for the shape of a random translation surface.
The goal of this paper is to study the asymptotic growth rate of the 
expected value of the covering radius of a translation surface as a function of genus.  
This is the maximal distance of any point to a singularity. In the case of the minimal 
stratum the covering radius is the same as the diameter, up to a factor $2$.   
A motivation for this study is a paper of Mirzakhani \cite{Mirzakhani} in which she computed
the expected value of several geometric functions (such as systole, Cheeger constant, etc.) 
on $\cM_g$, the moduli space
of Riemann surfaces of genus $g$, equipped with the Weil--Petersson volume measure
$\nu_{\mathrm{wp}}$. For example, she proved that the expected value of the diameter of a generic hyperbolic 
surface of genus $g$ grows like $\log g$ as $g\to \infty$. Specifically
\begin{equation*}
\mathbb{E}_{\cM_g} (\diam )   
  = \frac{\displaystyle{\int_{\cM_g} \diam  \, d\nu_{\mathrm{wp}} }} {\Vol_{\mathrm{wp}} ( \cM_g )} 
  \asymp \log g
\end{equation*}
where $\asymp$ means that the two sides are equal up to uniform multiplicative constants 
that are independent of $g$. 

The space of all translation surfaces is naturally stratified by the number and the type
of the singularities they can have and the expected shape of a translation surface may be 
different depending on the stratum.

To get some possible idea about shapes, a first guess
would be to translate the result of Mirzakhani directly. Namely, a hyperbolic
surface $x$ of genus $g$ has an area comparable to $g$. To make $x$ have area $1$, 
one needs to scale~$x$ down by a factor comparable to $\frac 1{\sqrt g}$. Then, the result 
of Mirzakhani would suggest that the expected value of the diameter should be comparable 
to $\frac{\log g }{\sqrt g}$.
However, the answer we find is different from this expected value.

Let $\kappa = (k_1, \dots, k_\ell)$ be a tuple of positive integers and let $\cH_1(\kappa)$
be the stratum of unit area translation surfaces with $\ell$ singularities of degrees
$k_1, \dots, k_\ell$.
Let $\nu$ be the normalized Lebesgue measure on a stratum $\cH_1(\kappa)$ as in
\cite{Masur-82,Veech-82}.

For a translation surface $(X,\omega) \in \cH(\kappa)$, the maximum distance from a point in $X$ to the set of singularities of $(X, \omega)$ is called the \emph{covering radius} of $(X, \omega)$.
This is equal to the maximum radius of an immersed Euclidean disk in $(X, \omega)$.
We denote the covering radius of $(X, \omega)$ by~$\crad(X, \omega)$.

\begin{introthm}[Expected covering radius] \label{introthm:covering_radius} 
Let $\cH(\kappa)$ be a stratum of translation surfaces of genus $g$. 
Then, for large values of $g$, we~have
 \begin{equation*}
  \mathbb{E}_{\cH_1(\kappa)} (  \crad) 
  = \frac{\displaystyle \int_{\cH_1(\kappa)} \crad(X) \, d\nu(X)}{\nu \big(\cH_1(\kappa)\big)}
  \leq 20 \cdot \sqrt{\frac{\log g}{g}}
  . 
 \end{equation*} 
\end{introthm}

In the special case of the minimal stratum $\cH_1(2g-2)$ where translation surfaces have exactly one singularity, the covering radius is up to a factor of $2$ the same as the diameter. We have therefore the following theorem which is a direct corollary of Theorem~\ref{introthm:covering_radius}.

\begin{introthm}[Expected diameter] \label{introthm:diameter_goes_to_zero} 
 For large values of $g$, we have
 \begin{equation*}
  \mathbb{E}_{\cH_1(2g-2)} (  \diam) 
  = \frac{\displaystyle \int_{\cH_1(2g-2)} \diam(X) \, d\nu(X)}{\nu \big(\cH_1(2g-2)\big)}
  \leq 40 \cdot \sqrt{\frac{\log g}{g}}
  .
 \end{equation*} 
\end{introthm}

This shows in particular that the expected value of the diameter in $\cH_1(2g-2)$ is smaller than what you would get from scaling a hyperbolic surface (by a factor $\frac 1 {\sqrt g}$) to have area $1$. So, the situation is different from that of hyperbolic surfaces.

It would be interesting to compute the expected value of the diameter of surfaces
in $\cH_1(\kappa)$ but that can not be achieved with our current methods. 
Unlike for the covering radius, it is not even clear at the moment whether or not the expected value of the diameter in every stratum goes to zero as the genus of the underlying surface 
goes to infinity. 

In contrast with \autoref{introthm:covering_radius} and \autoref{introthm:diameter_goes_to_zero},
we have the following absolute lower bound for the 
covering radius of any elements in $\cH_1(\kappa)$. 

\begin{introprop}[Lower bound on diameter]
For every $(X,\omega) \in \cH_1(\kappa)$, we have
\[
\crad(X) \geq \sqrt{\frac{2}{3\sqrt{3} \cdot (2g+\ell-2)}}.
\]
 
\begin{proof}
Let $\kappa = (k_1, \dots, k_\ell)$ and let $(X, \omega) \in \cH_1(\kappa)$. 
 Consider a Delaunay triangulation of~$(X, \omega)$. The number of triangles is 
 $2(2g+\ell-2)$. Hence, the area of the largest triangle is greater than or equal to 
 $\frac{1}{2(2g+\ell-2)}$. This triangle is inscribed in a circle of radius at least 
 $\sqrt{\frac{2}{3\sqrt{3} \cdot (2g+\ell-2)}}$.
 As the triangulation is Delaunay, the corresponding disk is an immersed Euclidean 
 disk in $X$. Hence its radius is a lower bound for the covering radius of the 
 translation surface. 
\end{proof}
\end{introprop}

\begin{rem}[Non-connectedness of $\cH_1(\kappa)$]
The stratum $\cH_1(\kappa)$ is not always connected. 
For $\kappa=(k_1, \dots, k_\ell)$, when every $k_i$ is even, there are different 
components corresponding to even and odd spin structure. Also, $\cH_1(2g-2)$ and 
$\cH_1(g-1, g-1)$ have a component consisting of hyperelliptic surfaces. 
So, a stratum may have up to three connected components 
(see \cite[Theorem 1]{kontsevich_zorich_03} for an exact statement). 
In our formula for the expected value of the diameter or the 
covering radius, we will not consider each component separately.
In fact, it is known that the volume of the hyperelliptic component
of $\cH_1(2g-2)$ is approximately of order $(2g)^{-2g}$ \cite[Theorem 1.1]{athreya_eskin_zorich_16}.
So we cannot even conclude that the 
expected value of the diameter for the hyperelliptic component goes to zero as $g\to \infty$. 
It was recently shown that the components associated to  even and odd spin structures 
have asymptotically equal volumes as $g\to \infty$ \cite{chen_moeller_sauvaget_zagier_19}. 
Hence, the statement of \autoref{introthm:diameter_goes_to_zero} does 
hold for these two components with a different constant. 
\end{rem} 

We now outline the proof of \autoref{introthm:covering_radius} which immediately implies 
\autoref{introthm:diameter_goes_to_zero}. For every $(X, \omega) \in \cH_1(\kappa)$, 
we find either an embedded disk or a cylinder that approximates the covering radius. When 
there exists a large embedded disk, we take out a parallelogram whose area is proportional 
to the area of the disk and glue the opposite sides to build a new translation surface.
The resulting translation surface is in the stratum $\cH(\kappa, 2)$
and its area is smaller than $(X, \omega)$ by a definite amount. 
We then renormalize this translation surface to have unit area.
We call this process taxing.  Because of the renormalization process,
the Jacobian of the taxing map is very large but the volume of $\cH(\kappa)$
and $\cH(\kappa, 2)$ are comparable. Hence, the volume of the subset of 
 $\cH(\kappa)$ where there is a large embedded disk is small. 
This allows us to show that the integral of the covering radius on these 
sets is small.

When there is a cylinder of large height, then either the cylinder has large area 
or a small circumference. In these cases, we bound the measure of the set of translation
surfaces that have such cylinders by bounding the associated Siegel--Veech constant.
Namely, for given length $\delta>0$, area $A \in [0, 1)$ and a translation surface 
$(X, \omega)\in \cH_1(\kappa)$, let $N_{\mathrm{cyl}}(X, \delta, A)$ be the number of cylinders 
in $X$ where the circumference is at most $\delta$ and the area is at 
least~$A$.

\begin{introthm}[Expected number of cylinders] \label{Thm:cyl}
There exists a constant $C>0$ such that for large values of $g$, we have 
\[
\mathbb{E}_{\cH_1(\kappa)} \big( N_{\mathrm{cyl}}(\param, \delta, A) \big)= 
\frac{\displaystyle \int_{\cH_1(\kappa)} N_{\mathrm{cyl}}(X, \delta, A) \, d\nu (X)} 
{\nu\big( \cH_1(\kappa) \big)} 
  \leq C \cdot g \cdot \delta^2 \cdot (1-A)^{2g+\ell-3}. 
\]
In particular, for $\Htc(\delta, A)$ the set of translation surfaces $(X,\omega) \in \cH_1(\kappa)$ for which $N_{\mathrm{cyl}}(X, \delta, A)$ is not zero, we have
\begin{equation*}
 \frac{\nu\left(\Htc(\delta, A)\right)}{\nu\left(\cH_1(\kappa)\right)}
 \leq C \cdot g \cdot \delta^2 \cdot (1-A)^{2g+\ell-3}
 .
\end{equation*}
\end{introthm}

For $A=0$, this is analogous to the estimate given by Mirzakhani 
\cite[Theorem 4.2]{Mirzakhani} for the Weil--Petersson volume of the set $\cM_g^\delta$ of Riemann surfaces of genus $g$ with 
at least one closed curve of length less than or equal to $\delta$. Namely
\[
\frac{\Vol_{\mathrm{wp}} (\cM_g^\delta)}{\Vol_{\mathrm{wp}} (\cM_g)} \asymp \delta^2. 
\]
In the setting of translation surfaces, another notion of thin part is the set of translation 
surfaces that have a short saddle connection.
In fact, for a stratum $\cH(\kappa)$ of translation surfaces and $\Ht(\delta)$ the set 
of translation surfaces in $\cH(\kappa)$ that have a saddle connection of length at 
most $\delta$, Masur and Smillie showed (compare equation (7) in the proof of 
Theorem~10.3 in \cite{masur_smillie_91})
\[
\nu\big( \Ht(\delta) \big) = O(\delta^2)
.
\]
However, the dependence of the constant on the genus or more generally on the stratum 
was not known. For a complete treatment of this topic, we also find an estimate for the number 
$N_{\mathrm{sc}}(X,\delta)$ of saddle connections of length at most $\delta$ in $X$. 

\begin{introthm}[Expected number of saddle connections] \label{introthm:sc}
There exists a constant $C'>0$ such that for large values of $g$, we have
\[
\mathbb{E}_{\cH_1(\kappa)} \big( N_{\mathrm{sc}}(\param, \delta) \big)= 
\frac{\displaystyle \int_{\cH_1(\kappa)} N_{\mathrm{sc}}(X, \delta) \, d\nu(X)}
  {\nu\big( \cH_1(\kappa) \big)}\leq C' \cdot g^2 \cdot \delta^2. 
\]
In particular, for $\Ht(\delta)$ the set of translation surfaces $(X,\omega) \in \cH_1(\kappa)$
for which $N_{\mathrm{sc}}(X, \delta)$ is not zero, we have
\begin{equation*}
\frac{ \nu\left(\Ht(\delta)\right)}{\nu\big(\cH_1(\kappa)\big)}
 \leq C' \cdot g^2 \cdot \delta^2.
\end{equation*}
\end{introthm}

\begin{rem}[Siegel--Veech constants and non-connectedness]
Theorems \ref{Thm:cyl} and \ref{introthm:sc} are proven in the appendix and in the proof, 
we make use of Siegel--Veech constants. Explicit formulas for 
values of various Siegel--Veech constants were computed by Eskin--Masur--Zorich in 
\cite{eskin_masur_zorich_03} in terms of combinatorial data and volumes of related 
strata with lower complexity. However, their methods give a precise
answer only for connected strata. Hence, we do not compute exact values, rather we find 
upper bounds for Siegel--Veech constants that suffice for our purposes.
Recently, Aggarwal in \cite{aggarwal_18} used the recursive formula for volume of 
strata given by Eskin--Okounkov \cite{EO} to compute the asymptotic growth rate of these 
volumes (see also \cite{chen_moeller_zagier_18} for the principal stratum and 
\cite{sauvaget_18} for the minimal~stratum). 

One should also compare \autoref{Thm:cyl} and \autoref{introthm:sc} to computations 
for values of various Siegel--Veech constants given in the appendix in 
\cite{aggarwal_18} written by Anton Zorich. For example, 
\autoref{Thm:cyl} is very similar to \cite[Corollary 5]{aggarwal_18}. The difference is that
in \cite[Corollary~5]{aggarwal_18}, only saddle connections bounding a cylinder of 
multiplicity $1$ are counted. However, higher multiplicity saddle connections 
do not pose a problem; this has been made precise by Aggarwal in \cite{aggarwal_18_SV}.
Also, the assumption on area being at least~$A$ contributes a factor of $(1-A)^{2g-2}$ 
to the estimate which is consistent with a result of Vorobets 
\cite[Theorem 1.8]{vorobets_2005}. Hence, \autoref{Thm:cyl} and \autoref{introthm:sc}
essentially follow from a combination of these results. However, we write a details
proof in the case of $\cH_1(2g-2)$, namely, we show that by a careful reading of 
\cite{eskin_masur_zorich_03} and incorporating the estimates 
given in \cite{aggarwal_18} one can obtain these theorems. 
\end{rem}

\paragraph{Acknowledgements}
We would like to thank Jon Chaika for suggesting that the diameter of a generic 
translation surface of high genus may go to zero which sparked this project. 
We are indebted to Anton Zorich for very helpful discussions on calculations of
Siegel--Veech constants. We also thank Jon Chaika, Samuel Lelièvre, 
Martin M\"oller, and Anton Zorich for their interest, helpful discussions, and pointing out 
several references. Furthermore, we want to thank the referees for their careful reading and helpful suggestions on improving the paper. The first author ac\-knowledges support from NSF grant DMS 1607512, the second author from NSERC Discovery grant RGPIN 06486, and the third author from NSERC grant RGPIN 06521. Some of this work took place at the Fields Institute during the Thematic Program on Teichm\"uller Theory and its Connections to 
Geometry, Topology and Dynamics.

\section{Four types of translation surfaces in \texorpdfstring{$\cH_1(\kappa)$}{a stratum}}

For this paper, a \emph{translation surface} $(X,\omega)$ is defined by a compact 
connected Riemann surface~$X$, a finite set $\Sigma \subseteq X$, and a translation 
structure on $X \setminus \Sigma$, i.e.\ a maximal atlas on $X \setminus \Sigma$ such 
that the transition maps are translations. The second parameter $\omega$ refers to the 
unique Abelian differential that is associated to a given translation structure on 
$X \setminus \Sigma$.
The elements of $\Sigma$ correspond to zeros of $\omega$ and 
are called \emph{singularities} of $(X,\omega)$. Every singularity $\sigma \in \Sigma$ is 
a cone point of the translation structure with \emph{cone angle} $2\pi(k+1)$ where $k$ 
is the order of $\sigma$ as a zero of $\omega$. The total sum of the orders of the zeros 
is equal to $2g-2$ where $g$ is the genus of~$X$. The translation structure defines also
a metric $d$ on $X$. With this metric, the \emph{diameter} of $X$ is defined to be 
$\diam(X) \coloneqq \max\{d(x,y): x, y\in X\}$
and the \emph{covering radius} of $X$ is defined to be $\crad \coloneqq \max \{d(x,\sigma): x\in X, \sigma\in \Sigma\}$.
See \cite{strebel_84} and \cite{zorich_06} for background
information on translation surfaces. 

Given a partition of $2g-2$ as a sum of integers $k_i\geq 1$, the \emph{stratum} 
$\cH(k_1,\ldots,k_\ell)$ is defined to be the set of all translation surfaces $(X,\omega)$ 
of genus $g$ with $\ell$ singularities of orders $k_1,\ldots,k_\ell$. The subsets of translation 
surfaces of area $1$ and area at most $1$ in $\cH(k_1,\ldots,k_\ell)$ are denoted by 
$\cH_1(k_1,\ldots,k_\ell)$ and $\cH_{\leq 1}(k_1,\ldots,k_\ell)$, respectively. 

A \emph{saddle connection} of $(X,\omega)$ is a geodesic segment from one 
singularity of $(X, \omega)$ to another (not necessarily different) 
singularity that is disjoint from $\Sigma$ in its interior. The translation structure 
of~$(X, \omega)$ associates a vector in $\mathbb{C}$, called the \emph{holonomy vector}, to a given oriented saddle connection. 
The holonomy vector of a saddle connection $\gamma$ can be expressed as  $\int_\gamma \omega$.
A \emph{cylinder} in $(X,\omega)$ of \emph{circumference} $k>0$ and \emph{height} $a>0$ is an isometrically embedded Euclidean cylinder $\mathbb{R}/k\mathbb{Z} \times (0,a)$. We assume that every cylinder is maximal. This implies that there are saddle connections on both boundary components of the cylinder.

A saddle connection can also be thought of as an element of the relative homology group
of~$(X, \omega)$ relative to $\Sigma$. If $\cB$ is a set of saddle connections that form a basis
for the relative homology, then the holonomy vectors of the elements of
$\cB$ determine $(X, \omega)$. For every $(X,\omega)$ and $\cB$, there is a neighborhood $U$ of 
$(X, \omega)$ in $\cH(k_1,\ldots,k_\ell)$ such that for every $(X', \omega')$ in $U$,
all elements of~$\cB$ (thought of as elements in the relative homology group) can still
be represented in $(X', \omega')$ as saddle connections. Then the set of holonomy vectors
of saddle connections in $\cB$ give coordinates for points in $U$. 
We refer to this set of holonomy vectors as \emph{period coordinates} for $\cH(k_1,\ldots,k_\ell)$ around $(X,\omega)$.
(See \cite{Masur-82} for details.) 

The period coordinates give an embedding from $U$ to 
$\mathbb{C}^{2g+\ell-1}$. The associated pullback measure in $\cH(k_1,\ldots,k_\ell)$
is called the \emph{normalized Lebesgue measure} denoted by $\nu$ and was studied by Masur \cite{Masur-82} and Veech \cite{Veech-82}. This also defines
a measure on $\cH_1(k_1,\ldots,k_\ell)$ in the following way. For an open set 
$V \subseteq \cH_1(k_1,\ldots,k_\ell)$, the measure is defined to be $\nu(U)$ where $U$ is the cone over $V$
of translation surfaces of area at most $1$. We abuse notation and denote the 
measure on $\cH_1(k_1,\ldots,k_\ell)$ also by $\nu$. 

In the following, let $\kappa= (k_1,\ldots, k_\ell)$ be a fixed partition of $2g-2$ and $\cH(\kappa)$ the corresponding stratum.
We show the main theorem by dividing the stratum $\cH_1(\kappa)$ into four (not necessarily disjoint) parts and investigate the behaviour of the expected value of the covering radius separately for every part.

\begin{definition}[$\Hsd$, $\Hpc$, $\Hrc$, and $\Hrd$] \label{def:subsets_of_stratum}
 For $g\geq 2$, define the following four subsets of $\cH_1(\kappa)$:
 \begin{itemize}
  \item Let $\Hsd$ be the subset of points $(X, \omega)$ where we have $\crad(X) < 18 \cdot \sqrt{\frac{\log g}{g}}$.
  \item Let $\Hpc$ be the subset of points $(X, \omega)$ where $\crad(X) \geq \frac{1}{\sqrt{g}}$ and where there exists a cylinder $C(X)$ of height at least $\crad(X)$ and $\area(C(X)) \leq \frac{1}{g}$. We call this the poor cylinder case.
  \item Let $\Hrc$ be the subset of points $(X, \omega)$ where $\crad(X) \geq \frac{1}{\sqrt{g}}$ and where there exists a cylinder $C(X)$ of height at least $\crad(X)$ and $\area(C(X)) \geq \frac{1}{g}$. We call this the rich cylinder case.
  \item Let $\Hrd$ be the subset of points $(X, \omega)$ where $\crad(X) \geq 18 \cdot \sqrt{\frac{\log g}{g}}$ and where there exists an embedded disk $D(X)$ of diameter at least $\crad(X)$. We call this the rich disk case.
 \end{itemize}
\end{definition}

Note that for a given translation surface $(X,\omega)$, we can be in the cylinder case and in the rich disk case.
Moreover, the choice of $C(X)$ in the cylinder case and $D(X)$ in the rich disk case is not canonical. However, for every translation surface $(X,\omega)$ in one of $\Hpc$, $\Hrc$, or $\Hrd$, we fix $C(X)$ or $D(X)$, respectively.
In particular, we fix $D(X)$ such that the holonomy vector defining the location of the center of $D(X)$ is locally constant in~$\Hrd$ (see the proof of \autoref{lem:richdisk_taxing_isometric_embedding} for details on this choice).

\begin{lem}[The four cases cover $\cH_1(\kappa)$]
 For $g\geq 2$, we have
 \begin{equation*}
  \cH_1(\kappa) = \Hsd \cup \Hpc \cup \Hrc \cup \Hrd
  .
 \end{equation*}
 
\begin{proof}
 Let $(X, \omega) \in \cH_1(\kappa) \setminus \Hsd$. Then there exists a point $x \in X$ such that $d(x,\Sigma) = \crad(X) \geq 18 \cdot \sqrt{\frac{\log g}{g}}$. In particular, there exists an immersed, locally flat, open disk around this point with radius $\crad(X)$.
 
 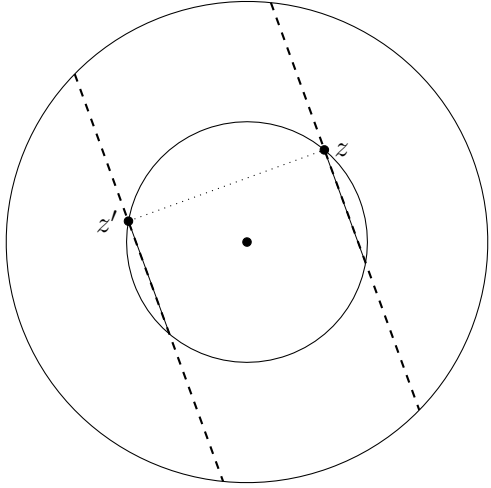
\begin{figure}
 \begin{center}
  \begin{tikzpicture}[scale=0.8, rotate=110]
   \draw[fill] (0,0) circle (2pt);
   \draw (0,0) circle (4cm);
   \draw (0,0) circle (2cm);
   \draw[fill] (60:2) circle (2pt) node[left]{$z'$};
   \draw[fill] (-60:2) circle (2pt) node[right]{$z$};
   \draw[dotted] (60:2) -- (-60:2);
   
   \draw[thick, dashed] (25.7:4) -- (154.3:4);
   \draw[very thin] (60:2) -- (120:2);
   \draw[thick, dashed] (-25.7:4) -- (-154.3:4);
   \draw[very thin] (-60:2) -- (-120:2);
  \end{tikzpicture}
  \caption{If $z$ and $z'$ are identified in the smaller disk, then the dotted line between them is a core curve of a cylinder. Hence, the two dashed lines are identified and give a lower bound on the height of the cylinder.}
  \label{fig:disk_wrapped_around_cylinder}
 \end{center}
 \end{figure}
 
 Consider an immersed disk with the same center but radius $\frac{1}{2} \crad(X)$. If this disk is not embedded, then two points in the disk have to be identified. This defines a closed geodesic and hence a core curve of a cylinder. The circumference of this cylinder is the length of the core curve which is at most $\crad(X)$. The height of this cylinder then has to be at least~$\crad(X,\omega)$ (see \autoref{fig:disk_wrapped_around_cylinder}). Furthermore, we have $\crad(X) \geq 18 \cdot \sqrt{\frac{\log g}{g}} \geq \frac{1}{\sqrt{g}}$.
 
 So, if we are not in the rich disk case then we are in the (poor or rich) cylinder case. 
\end{proof}
\end{lem}

\section{The poor cylinder case}

We will deal with the poor cylinder case and the rich cylinder case similarly by 
estimating the measure of the subset of $\Hpc$ (or $\Hrc$) where the height of the cylinder
is in a certain given range. Essentially, we use a Riemann sum argument 
to estimate the integral. 
We will need the following result proven in the appendix.

\begin{repcor}{lem:measure_cylinders_small_circumference}[Measure of $\Htc(\delta)$]
 There exists a constant $C>0$ such that for large values of $g$ and
 $\Htc(\delta)$ the set of translation surfaces $(X,\omega) \in \cH_1(\kappa)$ for which $N_{\mathrm{cyl}}(X, \delta, 0)$ is not zero, we have
 \begin{equation*}
  \frac{\nu\left(\Htc(\delta)\right)}{\nu\left(\cH_1(\kappa)\right)}
   \leq C \cdot g \cdot \delta^2 .
 \end{equation*}
\end{repcor}

Recall that for $(X, \omega) \in \Hpc$, the height of $C(X)$ is at least 
$\crad(X) \geq \frac{1}{\sqrt{g}}$ and the area of~$C(X)$ is at most 
$\frac{1}{g}$. We divide $\Hpc$ into subsets based on the height of the cylinder. 
For every~$n\geq 1$, define
\begin{equation*}
 \Hpc(n) =
 \left\{ (X,\omega) \in \Hpc : \height(C(X)) \in
   \left(2^{n-1} \cdot \frac{1}{\sqrt{g}} , 2^n \cdot \frac{1}{\sqrt{g}} \right) \right\}
 .
\end{equation*}
In particular, if $(X, \omega)$ is a translation surface in $\Hpc(n)$ then we have that $(X,\omega)$ contains a cylinder whose circumference is at most $\frac{1}{2^{n-1} \cdot \sqrt{g}}$.
With the above corollary, we can calculate the measure of $\Hpc(n)$.
\begin{cor}[Measure of $\Hpc(n)$] \label{cor:measure_thincylinder_n}
 For large values of $g$ and the constant $C$ from \autoref{lem:measure_cylinders_small_circumference}, we have
 \begin{equation*}
  \nu\left(\Hpc(n) \right)
  \leq \nu\left(\Htc\left( \frac{1}{2^{n-1} \cdot \sqrt{g}} \right)\right)
  \leq \frac{C}{2^{2n-2}} \cdot\nu\left(\cH_1(\kappa)\right)
  .
 \end{equation*}
\end{cor}

\begin{thm}[Expected covering radius on $\Hpc$] \label{thm:expected_diameter_Hpc}
 For large values of $g$ and the constant $C$ from \autoref{lem:measure_cylinders_small_circumference}, we have
 \begin{equation*}
  \frac{ \int_{\Hpc} \crad \, d\nu}{\nu\left(\cH_1(\kappa)\right)}
   \leq 4C \cdot \frac{1}{\sqrt{g}}
  .
 \end{equation*}
 
\begin{proof}
 A cylinder contains a singularity on each of its two boundary components. The distance from a point half way across the cylinder to the set of singularities is at least one half of the height of the cylinder, so the height is at most twice the covering radius.
 Hence, for an $(X,\omega) \in \Hpc(n)$, we have
 $\crad(X) \in \left(2^{n-2} \cdot \frac{1}{\sqrt{g}}, \, 2^{n} \cdot \frac{1}{\sqrt{g}} \right)$.  
 We can now calculate the integral by using $\Hpc = \cup_{n\geq 1} \Hpc(n)$.
 For the third inequality below we use \autoref{cor:measure_thincylinder_n}.
 \begin{align*}
  \int_{\Hpc} \crad \, d\nu 
  &\leq \sum_{n=1}^\infty \int_{\Hpc(n)} \crad \, d\nu \\
  &\leq \sum_{n=1}^\infty \nu\left(\Hpc(n)\right) \cdot 2^{n} \cdot \frac{1}{\sqrt{g}} \\
  &\leq \sum_{n=1}^\infty \frac{C}{2^{2n-2}} \cdot \nu\left(\cH_1(\kappa)\right) \cdot 2^{n} \cdot \frac{1}{\sqrt{g}} \\
  &=   4C \cdot \frac{1}{\sqrt{g}} \cdot \nu\left(\cH_1(\kappa)\right) \cdot \sum_{n=1}^\infty \frac{1}{2^{n}} \\
  &=   4C \cdot \frac{1}{\sqrt{g}} \cdot \nu\left(\cH_1(\kappa)\right)
 \end{align*}
 
 This finishes the proof of the statement.
\end{proof}
\end{thm}

\section{The rich cylinder case}

We can do a similar approach for the rich cylinder case as in the poor cylinder case.
Recall that for $(X, \omega) \in \Hrc$, the height of the cylinder $C(X)$ is at least 
$\crad(X) \geq \frac{1}{\sqrt{g}}$ 	
and the area of $C(X)$ is at least $\frac{1}{g}$. 

For $n, m \geq 1$, consider $\Hrc(n,m)$ to be the subset of $\Hrc$ where the height of $C(X)$ 
is in $\left(2^{n-1} \cdot \frac{1}{\sqrt{g}}, 2^n \cdot \frac{1}{\sqrt{g}}\, \right)$ and the area is in $\left( m \cdot \frac{1}{g}, (m+1) \cdot \frac{1}{g}\, \right)$.
The area of a translation surface in $\Hrc(n,m)$ is $1$, therefore $\Hrc(n,m)=\emptyset$ for $m \geq g$. We will state bounds on the volume of $\Hrc(n,m)$ for all $m\geq 1$ noting that the statements are trivially true for $m \geq g$.

Since the area of the cylinder $C(X)$ for an $(X,\omega) \in \Hrc(n,m)$ is bounded from above by $(m+1) \cdot \frac{1}{g}$ and the height is bounded from below by $2^{n-1} \cdot \frac{1}{\sqrt{g}}$, the circumference of $C(X)$ is bounded from above by~$\frac{m+1}{2^{n-1}} \cdot \frac{1}{\sqrt{g}}$.

With \autoref{thm:measure_cylinders_large_area}, we can calculate the measure of $\Hrc(n,m)$.

\begin{cor}[Measure of $\Hrc(n,m)$] \label{cor:measure_thincylinder_n_m}
 For large values of $g$ and the constant $C$ from \autoref{thm:measure_cylinders_large_area}, and for all $n,m \geq 1$
 \begin{equation*}
  \nu\left(\Hrc(n,m) \right)
  \leq \frac{4C \cdot (m+1)^2}{2^{2n}} \cdot e^{-m} \cdot \nu(\cH_1(\kappa))
  .
 \end{equation*}
 
\begin{proof}
 This follows directly from \autoref{thm:measure_cylinders_large_area} with the following calculation:
\pagebreak[2]
\begin{align*}
 \nu( \Hrc(n, m) )
 & \leq \Htc \left( \frac{m+1}{2^{n-1}} \cdot \frac{1}{\sqrt{g}} , \, \frac{m}{g} \right) \\
 & \leq C \cdot g \cdot \frac{(m+1)^2}{2^{2n-2}} \cdot \frac{1}{g} \cdot \left( 1- \frac{m}{g} \right)^{2g+\ell-3} \cdot \nu(\cH_1(\kappa)) \\
 & \leq \frac{4C \cdot (m+1)^2}{2^{2n}} \cdot \left( e^{- \sfrac{m}{g}} \right)^{2g+\ell-3} \cdot \nu(\cH_1(\kappa)) \\
 & \leq \frac{4C \cdot (m+1)^2}{2^{2n}} \cdot e^{-m} \cdot \nu(\cH_1(\kappa)) \qedhere
\end{align*}
\end{proof}
\end{cor}

\begin{thm}[Expected covering radius on $\Hrc$] \label{thm:expected_diameter_Hrc}
 For large values of $g$ and the constant $C$ from \autoref{thm:measure_cylinders_large_area}, we have
 \begin{equation*}
  \frac{ \int_{\Hrc} \crad \, d\nu}{\nu\left(\cH_1(\kappa)\right)}
  \leq 44C \cdot \frac{1}{\sqrt{g}}
  .
 \end{equation*}
 
\begin{proof}
 Note that the height of $C(X)$ cannot be larger than twice the covering radius of $X$. Hence, for an $(X,\omega) \in \Hrc(n,m)$, we have
 $\crad(X) \in \left(2^{n-2} \cdot \frac{1}{\sqrt{g}}, \, 2^{n} \cdot \frac{1}{\sqrt{g}} \right)$.

 We can now calculate the integral by using $\Hrc = \cup_{m\geq 1} \cup_{n\geq 1} \Hrc(n,m)$.
 For the third inequality below we use \autoref{cor:measure_thincylinder_n_m}.
 \begin{align*}
  \int_{\Hrc} \crad \, d\nu 
  &\leq \sum_{m=1}^\infty \sum_{n=1}^\infty \int_{\Hrc(n,m)} \crad \, d\nu \\
  &\leq \sum_{m=1}^\infty \sum_{n=1}^\infty \nu\left(\Hrc(n,m)\right) \cdot 2^{n} \cdot \frac{1}{\sqrt{g}} \\
  &\leq \sum_{m=1}^\infty \sum_{n=1}^\infty \frac{4C \cdot (m+1)^2}{2^{2n}} \cdot e^{-m} \cdot \nu(\cH_1(\kappa)) \cdot 2^{n} \cdot \frac{1}{\sqrt{g}} \\
  &= 4C \cdot \frac{1}{\sqrt{g}} \cdot \nu\left(\cH_1(\kappa)\right) \cdot \sum_{m=1}^\infty  (m+1)^2 \cdot e^{-m} \cdot \sum_{n=1}^\infty \frac{1}{2^{n}} \\
  &= 4C \cdot \frac{1}{\sqrt{g}} \cdot \nu\left(\cH_1(\kappa)\right) \cdot \sum_{m=1}^\infty  (m+1)^2 \cdot e^{-m}
 \end{align*}
 
 We now find a bound for the sum in this term. Note that $( (m+1)^2 \cdot e^{-m} )_{m \geq 1}$ is a decreasing sequence. In particular, we have for every $m\geq 6$:
 \begin{equation*}
  \frac{(m+2)^2 \cdot e^{-(m+1)}}{(m+1)^2 \cdot e^{-m}}
  = \left( \frac{m+2}{m+1} \right)^2 \cdot e^{- 1}
  \leq \left( \frac{8}{7} \right)^2 \cdot \frac{1}{e}
  < \frac{1}{2}
  .
 \end{equation*}
 Hence, the sum is bounded by seven times the first term; that is, we have
 \begin{equation*}
  \sum_{m=1}^\infty  (m+1)^2 \cdot e^{-m}
  \leq 7 \cdot 2^2 \cdot  e^{- 1}
  < 11
  .
 \end{equation*}
 This finally shows the statement.
\end{proof}
\end{thm}

\section{The rich disk case}

Recall that $\Hrd$ is the subset of translation surfaces $(X, \omega)$ in $\cH_1(\kappa)$ where $\crad(X) \geq 18 \cdot \sqrt{\frac{\log g}{g}}$ and where there exists an embedded disk of diameter at least $\crad(X)$.

Note that the area of an embedded disk can never be larger than $1$, hence the diameter of the disk has to be smaller than $\frac{2}{\sqrt{\pi}}$. Hence on $\Hrd$, the covering radius of a translation surface is globally bounded by $\frac{1}{\sqrt{\pi}}$. In particular, for small genera we have $\frac{1}{\sqrt{\pi}} \leq 18 \cdot \sqrt{\frac{\log g}{g}}$ and hence $\Hrd$ is empty.

The idea in the rich disk case is to take some area from the large embedded disk by removing a parallelogram and to distribute this area to the rest of the surface. 
Note that when doing so, we leave the stratum and change the topology of the surface. We will show that we end up in a subset of the image stratum of small volume. We will argue that this forces the set of surfaces with large embedded disks to have small volume as well. In what follows, it will be convenient to consider surfaces of area less than or equal to $1$.  
Let $\Hrdc\subseteq \cH_{\leq 1}(\kappa)$ be the cone over the set $\Hrd \subseteq \cH_1(\kappa)$, that is, the set of translation surfaces
that are obtained from some $(X, \omega) \in \Hrd$ by scaling by a factor 
$\sqrt \lambda$ with $0 < \lambda \leq 1$. 

We now describe the parallelogram that is to be removed and explore its properties.
For the fixed $g\geq 2$, set $\xi = \sqrt{\frac{3}{2}} \cdot \sqrt{\frac{\log g}{g}}$.
Let $(X,\omega) \in \Hrd$ and let $(x_1,x_2,x_3) \in (-\xi,\xi)^6 \subseteq \mathbb{C}^3$.
By definition, there is a choice of an embedded disk $D(X)$ in $X$ with diameter
 \[
 d \geq \crad(X) \geq 18 \cdot \sqrt{\frac{\log g}{g}} 
 \] 
and center~$c$.
We consider the parallelogram $P = P(X, x_1, x_2, x_3)$ in $X$ with center $c' = c+x_1$ and edges $(3\xi, 0) + x_2$ and $(0, 3\xi) + x_3$.

\begin{lem}[Properties of $P(X, x_1, x_2, x_3)$] \label{lem:properties_removed_parallelogram}
 For $(X,\omega) \in \Hrd$ and $(x_1,x_2,x_3) \in (-\xi, \xi)^6$, the parallelogram $P=P(X, x_1, x_2, x_3)$ has the following properties:
 \begin{enumerate}
  \item The parallelogram $P$ is embedded.
  \item The edges of $P$ are shorter than any geodesic segment from a corner of $P$ to itself or to another corner (except possibly the diagonals of $P$).
  \item We have
   \begin{equation*}
    \frac{9}{2} \cdot \frac{\log g}{g} \leq \area(P(X, x_1, x_2, x_3)) \leq \frac{51}{2} \cdot \frac{\log g}{g}
    .
   \end{equation*}
 \end{enumerate}

\begin{proof}
 \begin{enumerate}
  \item First note that the distance between $c'$ and $c$ is at most $\sqrt{2} \xi$. Second, the distance between $c'$ and a corner of the parallelogram is half the length of the corresponding diagonal, that is, it is bounded above by $\frac{1}{2} \cdot \sqrt{(3\xi + \xi + \xi)^2 + (3\xi + \xi + \xi)^2} = \frac{5\sqrt{2}}{2} \xi$. Therefore, every corner of $P$ has distance at most
  \begin{equation*}
   \frac{7\sqrt{2}}{2}\xi
   = \frac{7\sqrt{3}}{2} \cdot \sqrt{\frac{\log g}{g}}
   < 9 \cdot \sqrt{\frac{\log g}{g}}
   \end{equation*}
   from $c$ and hence is contained in $D(X)$. In particular, the whole parallelogram $P$ is contained in the embedded disk and therefore is embedded itself.
  \item The length of an edge of $P$ is at most $\sqrt{(4\xi)^2 + \xi^2} = \sqrt{17}\xi$.
  Note also that all the geodesic segments between corners of $P$ within $D(X)$ are edges or diagonals of $P$.
  The distance from any corner of the parallelogram to the complement of $D(X)$ is at least
  \begin{equation*}
   \left( 9-\frac{7\sqrt{3}}{2} \right) \cdot \sqrt{\frac{\log g}{g}}
   .
  \end{equation*}
  Therefore, any loop from a corner to itself or to another corner that leaves the disk has length at least
  \begin{equation*}
   \left( 18-7\sqrt{3} \right) \cdot \sqrt{\frac{\log g}{g}}
   \geq \sqrt{\frac{51}{2}} \cdot \sqrt{\frac{\log g}{g}}
   = \sqrt{17}\xi
   ,
  \end{equation*}
  so it is longer than any edge of the parallelogram.
  \item The area of $P(X, x_1, x_2, x_3)$ can be calculated as the determinant of the matrix with entries $(3\xi, 0) + x_2$ and $(0, 3\xi) + x_3$. As $x_2, x_3 \in (-\xi, \xi)^2$, the determinant is bounded from above by $4\xi \cdot 4\xi - (-\xi) \cdot \xi = 17 \xi^2 = \frac{51}{2} \cdot \frac{\log g}{g}$ and bounded from below by $2\xi \cdot 2\xi - \xi \cdot \xi = 3 \xi^2 = \frac{9}{2} \cdot \frac{\log g}{g}$.
 \end{enumerate}
\end{proof}
\end{lem}

To define now the taxing map locally, fix $(X_0,\omega_0)\in \Hrdc$.
Let $c_0$ be the center of the embedded disk $D(X_0)$ in $(X_0, \omega_0)$.
Let $v$ be a vector connecting a singularity $\sigma\in \Sigma$ to $c_0$ that can be represented
as a geodesic segment.
We can choose the center $c$ of the rich disks in the translation surfaces in $\Hrd$ in such a way that the vector $v$ joining the point in $\Sigma$ to $c$ is constant in a neighborhood of~$(X_0,\omega_0)$ in $\cH(\kappa)$ (see the remark after \autoref{def:subsets_of_stratum}).

\begin{definition}[Taxing map]
Define
\[
 T \colon \Hrdc \times (-\xi, \xi)^6 \to \cH_{\leq 1}(\kappa,2)
\] 
in the following way.

For a fixed $(X_0,\omega_0) \in \Hrdc$, consider a neighborhood in $\Hrdc$. And for any $(X,\omega)$ in this neighborhood and a chosen $\lambda \in (0,1)$, consider $(X', \omega') \in \Hrdc$  with $\area(X') = \lambda$ such that $(X,\omega)$ is the scaled version of $(X',\omega')$ (we multiply every saddle connection in $(X', \omega')$ by a factor 
of $\frac 1 {\sqrt \lambda}$).
 
For an $(x_1,x_2,x_3) \in (-\xi,\xi)^6 \subseteq \mathbb{C}^3$, let $P'= \sqrt{\lambda} P$ be the image of $P=P(X,x_1,x_2,x_3)$ in $(X', \omega')$ under scaling by $\sqrt \lambda$.
Now, we define $(Y', \zeta')$ to be the translation surface where we remove $P'$ from $(X', \omega')$ and glue the two pairs of parallel boundaries of $P'$. This introduces a new singularity $\tau$ with cone angle $3 \cdot 2\pi$, hence $(Y', \zeta') \in \cH_{\leq 1}(\kappa,2) \coloneqq \cH_{\leq 1}(k_1,\ldots,k_\ell,2)$.
The glued parallel sides give rise to a pair of loops based at the singularity $\tau$ that are shortest under all loops based at $\tau$.
 
 As shown in \autoref{lem:properties_removed_parallelogram},
 \[
 \area(Y') \leq \left(1 - \frac{9}{2} \cdot \frac{\log g}{g} \right) \cdot \area(X'). 
 \]
 We rescale $(Y', \zeta')$ by this factor $\left(1 - \frac{9}{2} \cdot \frac{\log g}{g}\right)^{-\frac{1}{2}}$, and obtain a surface $T\big( (X', \omega'), (x_1,x_2,x_3)\big)$ that has area less than or equal to $1$.
\end{definition}

We study the map $T$ by considering it as the concatenation of the maps 
\[
T_P \colon \Hrdc \times (-\xi, \xi)^6 \to  \cH_{\leq 1}(\kappa,2)
\]
(removing $P'$) and 
\[
 T_S \colon \cH(\kappa,2) \to \cH(\kappa,2)
\] 
(scaling the surface by $\left(1 - \frac{9}{2} \cdot \frac{\log g}{g}\right)^{-\frac{1}{2}}$).

\begin{lem}[Properties of $T_P$] \label{lem:richdisk_taxing_isometric_embedding} \leavevmode
\begin{enumerate}
 \item Every translation surface in $\cH_{\leq 1}(\kappa,2)$ has at most $g$ preimages under $T_P$.
 \item We have $\Jac(T_P, ((X', \omega'), (x_1,x_2,x_3))) = (\area(X', \omega'))^3$ on $\Hrdc \times (-\xi,\xi)^6$.
\end{enumerate}

\begin{proof}
 Let $\ell$ be the number of singularities of $(X,\omega)\in \cH(\kappa)$. Then the complex dimension of $\cH_{\leq 1}(\kappa)$ is $2g+\ell-1$. Hence, $\Hrdc \times (-\xi,\xi)^6$ has complex dimension $2g+\ell+2$. Also, the complex dimension of $\cH_{\leq 1}(\kappa, 2)$ is 
 $2(g+1)+(\ell+1)-1 = 2g+\ell+2$. Hence the image and the domain have the same dimension.

We now show that a translation surface $(Y',\zeta')$ in $\cH_{\leq 1}(\kappa,2)$ has at most $g$ preimages under~$T_P$.
We first note that it may have no preimages at all. This might happen, for instance, if the area of~$(Y',\zeta')$ is close to $1$ and so the parallelogram that we will add later forces the area of the possible preimage to be larger than $1$. In that case $(Y',\zeta')$ is not in the image of $T_P$.  

Recall that $(Y', \zeta')$ has genus $g+1$, hence it can have at most $g$ singularities with cone angle $6\pi$. We now fix a singularity $\tau$ of $(Y', \zeta')$ with cone angle $6\pi$ and call the set of the other singularities~$\Sigma$. If the two shortest loops based at $\tau$ only intersect at $\tau$ then $(Y',\zeta')$ is possibly contained in the image of $T_P$. If the two shortest loops are not unique or intersect, $(Y',\zeta')$ cannot be in the image of $T_P$ for this given $\tau$. Thus assume the condition on two unique shortest loops at $\tau$ holds. We construct the only possible preimage $(X', \omega')$ of $(Y',\zeta')$ under~$T_P$ for the given singularity $\tau$. 

We cut open along these unique shortest loops at $\tau$ giving us the possibility to glue in a parallelogram $P'$. Doing so, all corners of the former parallelogram will be regular 
points (that is, their cone angle is $2\pi$) and we obtain a translation surface $(X', \omega')$.  The area may be greater than $1$. In this case again $(Y',\zeta')$ has no preimages.  Therefore assume $(X', \omega')$ has area~$\lambda \in(0,1]$, that is $(X', \omega') \in \cH_{\leq 1}(\kappa)$.
Let $(X, \omega) = \left(\frac 1{\sqrt \lambda} X', \frac 1{\sqrt \lambda} \omega'\right)$, that is $\area(X, \omega) = 1$. 

Let $c' \in (X,\omega) $ be the preimage of the point in $(X', \omega')$ which is the center of the parallelogram that was glued in to obtain~$(X', \omega')$. 

Let $y_2$ be the holonomy vector of the edge $e_2$ of the parallelogram $P'$ with the smaller imaginary part and $y_3$ be the holonomy vector of the edge $e_3$ of $P'$ with the larger imaginary part. 
Define $x_1 = c' - c$, $x_2 = \frac{1}{\sqrt \lambda} y_2 - (3\xi, 0)$, and $x_3 = \frac{1}{\sqrt{\lambda}} y_3 - (0, 3\xi)$.
Then $\big( \sqrt \lambda (X, \omega), (x_1, x_2, x_3)\big)$ is the only 
possible preimage of $(Y',\zeta')$ when $\tau$ is a fixed singularity with angle $6\pi$.

In particular, $T_P$ is locally injective. The equality of dimensions and the local injectivity imply that the image of $T_P$ is a subset of $\cH_{\leq 1}(\kappa,2)$ with non-empty interior.

To compare the measure in the domain and in the image, we locally (around the points
$(X',\omega')$ and $(Y', \zeta')$) choose \emph{compatible period coordinates} on 
$\Hrdc \subseteq \cH(\kappa)$ and on~$\cH(\kappa, 2)$.

For any triangulation $\Delta$ of $(Y',\zeta')$ by saddle connections, there is an open 
set $U_\Delta$ in~$\cH(\kappa,2)$ where the edges of $\Delta$ (as a topological triangulation) 
can still be represented by saddle connections. Also, any subset $\cE$ of edges of 
$\Delta$ that does not contain a contractible loop can be extended to a set $\cB$ of edges 
of $\Delta$ that form a basis for the relative homology of~$Y'$ relative to $\Sigma \cup \{\tau\}$.

Let $\cE$ be the set of the following three saddle connections in~$(Y', \zeta')$: 
the two shortest saddle connections $e_2$ and $e_3$ that connect $\tau$ to itself 
and a saddle connection~$e_1$ that connects a singularity in $\Sigma$ to $\tau$ and is disjoint from 
the other two such that the holonomy vector of $e_1$ is 
$y_1= v + x_1 - (x_2 + x_3)/2-(3\xi, 3\xi)/2$. (That is, $e_1$ connects a singularity in $\Sigma$ to the lower left corner of 
the parallelogram spanned by $e_2$ and $e_3$.)
 
Complete $\cE$ to a triangulation $\Delta$ of $(Y', \zeta')$ by successively adding edges
that are disjoint from previous edges. Let $\cE_\Sigma$ be the set of edges in $\Delta$
that start and end in singularities in $\Sigma$. Then $\cE \cup \cE_\Sigma$ spans the relative homology. 
This is because any edge $e$ connecting a singularity in $\Sigma$ to $\tau$ is either in $\cE$ or there is 
a contractible loop in $\Delta$ consisting of $e$, the edge $e_1$ in~$\cE$ connecting a singularity in 
$\Sigma$ to $\tau$, and an edge path in $\cE_\Sigma$. Let $\cB_\Sigma$ be a subset of 
$\cE_\Sigma$ so that $\cB = \cB_\Sigma \cup \cE$ forms a basis for the relative homology
of $Y'$ relative to $\Sigma \cup \{\tau\}$.

By construction, edges in $\cE_\Sigma$ can also be represented as saddle connections
in $(X', \omega')$. In fact, the edges associated to $\cB_\Sigma$ form a basis for
the relative homology of $X'$ relative to $\Sigma$ as~$|\cE| = 3$ and the complex 
dimension of $\cH_{\leq 1}(\kappa)$ is $3$ less than the dimension of $\cH_{\leq 1}(\kappa,2)$.
By making $U_\Delta$ smaller, we can ensure that any point in $\Hrdc$ that is a
preimage of a point in $U_\Delta$ is contained in an open set $U_X$ where edges in 
$\cB_\Sigma$ can still be represented by saddle connections. We call the period 
coordinates given by $\cB_\Sigma$ for points in $U_X$ and the period coordinates given 
by $\cB$ for points in $U_\Delta$ a pair of compatible period coordinates. 

On a pair of compatible period coordinates (together with the three vectors 
$x_1$, $x_2$, $x_3$), we therefore 
have that $T_P$ is affine: it is the identity on all but the last three vectors and
the coefficients of $e_1$, $e_2$ and $e_3$ depend only on $x_1$, $x_2$ and $x_3$. 
That is, $T_P$ can be represented in the following form 
\[ 
\begin{pmatrix}
I_{2g} & 0 \\
0 & A
\end{pmatrix} 
 \begin{pmatrix}
\cB_\Sigma\\
\cE
\end{pmatrix} 
+ \begin{pmatrix}
0\\
B
\end{pmatrix}.
\]
Recalling that
\begin{align*}
y_1&= \sqrt \lambda \left( x_1 + v -\bigg(\frac32 \xi, \frac32 \xi\bigg)-\frac{x_2}{2}-\frac{x_3}{2}\right)\\
y_2 &= \sqrt \lambda \big( x_2 + (3\xi, 0)\big)\\
y_3 &= \sqrt \lambda \big( x_3 + (0, 3\xi) \big)
\end{align*}
 we have  
\[
\begin{pmatrix}
y_1 \\
y_2\\
y_3 
\end{pmatrix} =
A \cE + B 
= \sqrt  \lambda \begin{pmatrix}
 1 &  -\frac{1}{2}& -\frac{1}{2}\\
0 & 1 & 0 \\
0 & 0 & 1
\end{pmatrix} 
\begin{pmatrix}
x_1 \\
x_2\\
x_3 
\end{pmatrix}
+ 
\sqrt  \lambda  \begin{pmatrix}
-(\frac32 \xi, \frac32 \xi) + v \\
(3\xi, 0)\\
(0, 3\xi)  
\end{pmatrix}
\]
The transformation above should be thought of as a map from 
$\mathbb{R}^6$ to $\mathbb{R}^6$. Hence, 
\begin{equation*}
\Jac(T_P, ((X',\omega'), (x_1,x_2,x_3))) = \sqrt{\lambda}^{\, 6} 
  = (\area (X', \omega'))^3.\qedhere
\end{equation*}
\end{proof}
\end{lem}

\begin{lem}[Jacobian of the scaling map] \label{lem:richdisk_taxing_Jacobian}
We have $\Jac(T_S) \geq g^{\sfrac{9}{2}}$ on $T_P(\Hrd \times (0,1])$.
 
\begin{proof}
 Let $(X, \omega) \in \cH_{\leq 1}(\kappa,2)$ be in the image of $T_P$. The map $T_S$ changes all of the $2 \cdot (2g+\ell+2) \geq 2 \cdot 2g$ real period coordinates uniformly by the factor $\left(1-\frac{\frac{9}{2} \cdot \log g}{g}\right)^{-\frac{1}{2}}$. For the calculation, we use that for every $x \in (0,1)$, we have $1 \geq 1-x^2$ and hence $\frac{1}{1-x} \geq 1 +x$.
 
 \begin{align*}
  \Jac(T_S) \ & = \left( \left( 1-\frac{\frac{9}{2} \cdot \log g}{g} \right)^{-\frac{1}{2}} \right)^{2 \cdot 2g} \\
  & = \left( \frac{1}{1- \frac{\frac{9}{2} \cdot \log g}{g} } \right)^{2g} \\
  & \geq \left( 1 + \frac{\frac{9}{2} \cdot \log g}{g} \right)^{2g} \\
  & \geq \left( 1 + \frac{1}{\frac{g}{\sfrac{9}{2} \cdot \log g}} \right)^{\left( \frac{g}{\sfrac{9}{2} \cdot \log g} + 1 \right) \cdot \frac{2g}{2\cdot \frac{g}{\sfrac{9}{2} \cdot \log g}}} \\
  & \geq e^{\sfrac{9}{2} \cdot \log g}  = g^{\sfrac{9}{2}}. \qedhere
 \end{align*}
\end{proof}
\end{lem}

Combining \autoref{lem:richdisk_taxing_isometric_embedding} and \autoref{lem:richdisk_taxing_Jacobian}, we get the following bound for the measure.

\begin{cor}[Measure of $\Hrd$]
 We have
 \begin{equation*}
  \nu \left( \Hrd \right)
  \leq \frac{2}{27} \cdot (\log g)^{-3} \cdot g^{-\sfrac{1}{2}} \cdot \nu\left(\cH_1(\kappa,2)\right)
  .
 \end{equation*}
 
\begin{proof}
 Recall that we use the same notation $\nu$ for the measure on the whole stratum $\cH(\kappa)$ and the measure on $\cH_1(\kappa)$.
 By definition, we have 
 \[
 \nu \left( \Hrd \right) = \nu \left( \Hrdc \right)
 \qquad\text{and}\qquad
 \nu\left(\cH_1(\kappa,2)\right) = \nu \left( \cH_{\leq 1}(\kappa, 2) \right).
 \]
 
 Denote the image of $\Hrdc \times (-\xi, \xi)^6$ under $T$ (or $T_P$) 
 with $\im(T)$ (or $\im(T_P)$, respectively). From \autoref{lem:richdisk_taxing_Jacobian} 
 and from the fact that $\im(T) \subseteq \cH_{\leq 1}(\kappa,2)$, we get
 \begin{equation*}
  \nu\left(\im(T_P)\right) \cdot g^{\sfrac{9}{2}}
  \leq \nu\left(\im(T)\right)
  \leq \nu\left( \cH_{\leq 1}(\kappa,2) \right)
  .
 \end{equation*}
 
 Let $\Hrdu$ be the set of translation surfaces in $\Hrdc$ that have area at least~$\frac{1}{2}$.
 As the complex dimension of $\Hrdc$ is $2g+\ell-1$, we have
 \begin{equation*}
  \nu\left( \Hrdu \right)
  = \left(1- \left( \frac{1}{2} \right)^{2(2g+\ell-1)} \right) \cdot \nu\left( \Hrdc \right)
  \geq \frac{1}{2} \cdot \nu\left( \Hrd \right)
  .
 \end{equation*}
 By \autoref{lem:richdisk_taxing_isometric_embedding}, we have
 \[
  g \cdot \nu\left(T_P\left(\Hrdu \times (-\xi, \xi)^6 \right) \right) \geq \left(\frac{1}{2}\right)^3 \cdot \nu\left(\Hrdu \right) \cdot (2\xi)^6
 .
 \] 
 This implies
 \begin{align*}
  \nu \left(T_P \left(\Hrdc \times (-\xi, \xi)^6 \right) \right)
  & \geq  \nu \left(T_P\left(\Hrdu \times (-\xi, \xi)^6 \right)\right) \\
  & \geq  8 \xi^6 \cdot g^{-1} \cdot \nu\left(\Hrdu \right)  \\
  & \geq 4\xi^6 \cdot g^{-1} \cdot \nu\left( \Hrd \right) .
 \end{align*}

 Combining these measure comparisons and inserting $\xi = \sqrt{\frac{3}{2}} \cdot \sqrt{\frac{\log g}{g}}$, we can now deduce
 \begin{align*}
  \nu\left( \Hrd \right)
  & \leq \frac{1}{4} \xi^{-6} \cdot g \cdot \nu(T_P(\Hrdc \times (-\xi, \xi)^6 )) \\
  & \leq \frac{2}{27} \cdot \frac{g^3}{(\log g)^3} \cdot g \cdot g^{-\sfrac{9}{2}} \cdot \nu\left( \cH_{\leq 1}(\kappa,2) \right)
  . \qedhere
 \end{align*}
\end{proof}
\end{cor}

\begin{thm}[Expected covering radius on $\Hrd$] \label{thm:expected_diameter_Hrd}
 For large values of $g$, we have
 \begin{equation*}
  \frac{ \int_{\Hrd} \crad \, d\nu}{\nu\left(\cH_1(\kappa)\right)}
   \leq \frac{1}{45} \cdot \sqrt{\frac{\log g}{g}}
  .
 \end{equation*}
 
\begin{proof}
 Recall that on $\Hrd$, the covering radius of a translation surface is globally bounded by~$\frac{1}{\sqrt{\pi}}$.
 This gives us the following calculation.
 
 \begin{align*}
   \int_{\Hrd} \crad \, d\nu 
  & \leq \frac{1}{\sqrt{\pi}} \cdot \nu\left(\Hrd\right) \\
  & \leq \frac{1}{\sqrt{\pi}} \cdot \frac{2}{27} \cdot (\log g)^{-3} \cdot g^{-\sfrac{1}{2}} \cdot \nu\left(\cH_1(\kappa,2)\right) \\
  & \leq \frac{2}{27 \sqrt{\pi}} \cdot \sqrt{\frac{\log g}{g}} \cdot (\log g)^{-\frac{7}{2}} \cdot \nu\left(\cH_1(\kappa,2)\right)
 \end{align*}
 
 By \cite[Theorem 1.4]{aggarwal_18}, we have the estimate
 $\nu \left(\cH_1(\kappa) \right) = \frac{4}{ \prod_{i=1}^\ell (k_i +1)} \cdot (1+\cO(\frac{1}{g}))$
 for a stratum $\cH_1(\kappa)$ with $\kappa = (k_1,\ldots,k_\ell)$.
 In particular, this also gives us
 $\nu \left(\cH_1(\kappa,2) \right) = \frac{4}{3 \cdot \prod_{i=1}^\ell (k_i +1)} \cdot (1+\cO(\frac{1}{g}))$.
 Hence, for large values of $g$, we have
\[
\frac{\nu \left(\cH_1(\kappa,2)\right)}{\nu \left(\cH_1(\kappa)\right)} \leq \frac 12. 
\]
The calculation $\frac{1}{27 \sqrt{\pi}} \leq \frac{1}{45}$ finishes the proof.
\end{proof}
\end{thm}

\section{Proof of the main theorem}

Now we can put together the ingredients for the proof of our main theorem.

We have shown in the previous three sections that the expected value of the covering radius goes to zero on $\Hpc$, on $\Hrc$, and on $\Hrd$.
When looking at the explicit statements in \autoref{thm:expected_diameter_Hpc} and \autoref{thm:expected_diameter_Hrc}, it is also clear that the rates $4C \cdot \frac{1}{\sqrt{g}}$ and $44C \cdot \frac{1}{\sqrt{g}}$ for some constant $C$ are faster than $\frac{1}{2} \cdot \sqrt{\frac{\log g}{g}}$ for large values of $g$. In \autoref{thm:expected_diameter_Hrd}, the rate is $\frac{1}{45} \cdot \sqrt{\frac{\log g}{g}}$.

The only missing part is $\Hsd$. However, this set is defined so that the covering radius is smaller than $18 \cdot \sqrt{\frac{\log g}{g}}$, hence the expected value of the covering radius is also bounded by~$18 \cdot \sqrt{\frac{\log g}{g}}$.
Therefore, by summing up these four expected values, we have proven \autoref{introthm:covering_radius} that we state here again.

\begin{thm}[Expected covering radius]
 For large values of $g$, we have
 \begin{equation*}
  \mathbb{E}_{\cH_1(\kappa)} (  \crad ) 
  = \frac{\displaystyle \int_{\cH_1(\kappa)} \crad(X) \, d\nu(X)}{\nu \big(\cH_1(\kappa)\big)}
  \leq 20 \cdot \sqrt{\frac{\log g}{g}}
  .
 \end{equation*} 
\end{thm}

\section{Appendix}

In this appendix, we find an upper bound for the number of
cylinders on a translation surface where an upper bound on the circumference and a lower bound on the area of the cylinder are given. We also give an upper bound for the number of saddle connections where an upper bound on the length is given. 
We then use these bounds to give upper bounds for the measure of the thin part of the stratum~$\cH_1(\kappa)=\cH_1(k_1,\ldots,k_\ell)$.

Here, we understand the $\delta$--thin part in three different ways. 
First, we consider the set $\Htc(\delta)$ of translation surfaces of area $1$ on which there exists a cylinder that has circumference at most~$\delta$.
Second, we do a similar computation for the subset $\Htc(\delta, A)\subseteq \Htc(\delta)$ where we add the assumption that the area of the thin cylinder is bounded from below by some constant $A \in [0,\frac{1}{2}]$.
Third, for a complete treatment of the volume question, we consider the usual thin part, namely the set $\Ht(\delta)$ of translation surfaces that contain a saddle connection of length at most~$\delta$ that do not necessarily bound a cylinder.

Recall that the complex dimension $d$ of $\cH(\kappa)$ is equal to $2g+\ell-1$.

\begin{thm}[Expected number of cylinders] \label{thm:measure_cylinders_large_area}
 Consider a stratum $\cH(\kappa)$ of complex dimension $d$.
 Given $\delta>0$, $0<A<1$, and $(X, \omega)\in \cH_1(\kappa)$, let $N_{\mathrm{cyl}}(X, \delta, A)$ be the number of cylinders in $X$ where the circumference is at most~$\delta$ and the area is at least~$A$.
 There exists a constant $C>0$ such that for
 large values of $g$, we have
 \[
  \mathbb{E}_{\cH_1(\kappa)} \big( N(\param, \delta, A) \big)= 
  \frac{\int_{\cH_1(\kappa)} N_{\mathrm{cyl}}(X, \delta, A) \, d\nu (X)} 
  {\nu\big( \cH_1(\kappa) \big)} 
  \leq C \cdot g \cdot \delta^2 \cdot (1-A)^{d-2}. 
 \]
In particular, for $\Htc(\delta)$ the set of translation surfaces $(X,\omega) \in \cH_1(\kappa)$ for which $N_{\mathrm{cyl}}(X, \delta, A)$ is not zero, we have
\begin{equation*}
\frac{\nu\left(\Htc(\delta, A)\right)}{\nu\left(\cH_1(\kappa)\right)}
 \leq C \cdot g \cdot \delta^2 \cdot (1-A)^{d-2} 
 .
\end{equation*}
\end{thm}

From \autoref{thm:measure_cylinders_large_area}, we can directly deduce the measure of $\Htc(\delta) = \Htc(\delta, 0)$ where we do not have any restriction on the area of the cylinder. This is used in the proof of \autoref{cor:measure_thincylinder_n}.

\begin{cor}[Measure of $\Htc(\delta)$]\label{lem:measure_cylinders_small_circumference}
 For $C$ as in \autoref{thm:measure_cylinders_large_area} and $\Htc(\delta)$ the set of translation surfaces $(X,\omega) \in \cH_1(\kappa)$ for which $N_{\mathrm{cyl}}(X, \delta, 0)$ is not zero, we have
\begin{equation*}
\frac{\nu\left(\Htc(\delta)\right)}{\nu\left(\cH_1(\kappa)\right)}
 \leq C \cdot g \cdot \delta^2 .
\end{equation*}
\end{cor}

For the sake of completeness, we also give the corresponding statement on the number of saddle connections.

\begin{thm}[Expected number of saddle connections] \label{Thm:sc}
 For given $\delta>0$ and a translation surface $(X, \omega)\in \cH_1(\kappa)$, let $N_{\mathrm{sc}}(X, \delta)$ be the number of saddle connections in $X$ of length at most $\delta$.
 There exists a constant $C'>0$ such that for
 large values of~$g$, we have
 \[
 \mathbb{E}_{\cH_1(\kappa)} \big( N(\param, \delta) \big)= 
 \frac{\int_{\cH_1(\kappa)} N_{\mathrm{sc}}(X, \delta) \, d\nu(X)}
   {\nu\big( \cH_1(\kappa) \big)}\leq C' \cdot g^2 \cdot \delta^2. 
 \]
In particular, for $\Ht(\delta)$ the set of translation surfaces $(X,\omega) \in \cH_1(\kappa)$ for which $N_{\mathrm{sc}}(X, \delta)$ is not zero, we have 
\begin{equation*}
\frac{\nu\left(\Ht(\delta)\right)}{\nu\left(\cH_1(\kappa)\right)}
 \leq C' \cdot g^2 \cdot \delta^2 
 .
\end{equation*}
\end{thm}

The key tools for these computations is to compute upper bounds for the Siegel--Veech constants 
associated to the two different situations. Siegel--Veech constants were introduced to study asymptotics for various counting problems of saddle connections or cylinders. Such a counting problem could be the asymptotics as $L$ goes to $\infty$ of the number of saddle connections whose holonomy vectors have multiplicity $1$ and length at most $L$. Given a translation surface $(X, \omega)$ in $\cH_1(\kappa)$ and a counting problem we are interested in,  we can consider the set $V = V(X) \subset \mathbb{R}^2$ of  holonomy vectors of the collection of saddle connections of interest (see \autoref{subsec:C} for details).
There is an action of $\SL(2,\mathbb{R})$ on each (connected component of a) stratum as well as on~$\mathbb{R}^2$.   Then $V(X)$ is a set of vectors equivariant under these actions. The Siegel-Veech formula~\eqref{eq:transform} below holds for any $\SL(2,\mathbb{R})$--invariant measure; in particular for the normalized Lebesgue measure that we consider in this paper. 
Now, for a connected component $\cH$ of a stratum,
a compactly supported continuous function 
$f\colon\, \mathbb{R}^2\to \mathbb R$ (for example a characteristic function of a ball) 
and a set of vectors~$V(X)$ for every $(X, \omega) \in \cH$  define $\hat f\colon\, \cH\to \mathbb R$ by
\begin{equation} \label{Eq:f-hat} 
\hat f(X,\omega)=\sum_{v\in V(X)} f(v).
\end{equation} 
Then there is a constant $c(V)$ \cite[Theorem 0.5]{Veech_98}, called the \emph{Siegel--Veech 
constant}, which is independent of $f$ and such that 
\begin{equation}
 \label{eq:transform}
 \frac{1}{\nu(\cH)}\int_{\cH} \hat f \, d\nu
 = c(V) \int_{\mathbb {R}^2} f \, dxdy.
\end{equation}

We will not compute these Siegel--Veech constants but determine upper bounds 
which suffice to calculate upper bounds for the counting problems of interest.

For connected strata, Zorich in an appendix to \cite{aggarwal_18} gives estimates for the Siegel--Veech constants of cylinders without conditions on the area but with conditions on the bounding saddle connections, in particular that the multiplicity is $1$. Combining this with \cite{aggarwal_18_SV} where Aggarwal shows that the constants for higher multiplicity are negligible in comparison to multiplicity $1$, we have a first estimate for connected strata without the area constraint on cylinders. We now have to add arguments for the area of the cylinder, arguments for non-connected strata, sum up the different types of cylinders, and compare the bound for cylinders with the one for saddle connections.

For the sake of clarity, we use the first four subsections of the appendix to
carry out the computations in detail in the (non-connected) minimal stratum $\cH_1(2g-2)$. 
To do this, we follow the arguments in \cite{eskin_masur_zorich_03}, but consider additionally the area condition for cylinders. In particular, we show that the Siegel-Veech constant for cylinders with bounding saddle connections of multiplicity greater than $1$ is dominated by the one for multiplicity~$1$ as $g$ goes to $\infty$.

\subsection{Setting to calculate Siegel--Veech constants in the minimal stratum}

We now begin the computations in the case of the minimal stratum $\cH_1(2g-2)$.
The key idea of the proof from \cite{eskin_masur_zorich_03} is that a set of $p$ homologous short saddle connections gives a decomposition of the translation surface into several pieces that are themselves surfaces with boundary.
In our setting, there are three different types of surfaces with boundary 
than can appear in this construction (see \autoref{fig:types_subsurfaces}):

\begin{figure}
 \begin{center}
  \begin{tikzpicture}[scale=0.6]
   \newcommand\hole{
    \draw[very thick] (0,0) to[out=180,in=-90]
    (-2,1.5) to[out=90,in=180]
    (-1.2,2.3) to[out=0,in=100]
    (-0.1,1.3) to[out=-70,in=90,looseness=0.5] (0,0);
    \draw[fill] (0,0) circle (2pt);
   }
   
   \newcommand\genus{
    \draw[thick] ((1,0.05) to[bend left]  (2.5,0.05);
    \draw[thick] ((0.8,0.1) to[bend right] (2.7,0.05);
   }

   \hole
   \begin{scope}[rotate=180]
    \hole
   \end{scope}
   \draw[thick] (-2,1.6) to[out=-110,in=80]
   (-2.5,0.3) to[out=-100, in=90]
   (-5,-3.5) to[out=-90, in=180]
   (-3,-6) to[out=0, in=-90]
   (0,-4) to[out=90, in=-150]
   (1.6,-2.2);
   \begin{scope}[rotate=55, xshift=-5.5cm]
    \genus
   \end{scope}

   \begin{scope}[xshift=9cm]
    \hole
    \begin{scope}[rotate=180, xshift=-1.3cm, yshift=0.7cm]
     \hole
    \end{scope}
    \draw[thick] (0,0) .. controls +(-110:1.5cm) and +(-130:1.2cm) .. (1.3,-0.7);
    \draw[thick] (-2,1.6) to[out=-110,in=80]
     (-2.5,0) to[out=-100, in=80]
     (-5,-3) to[out=-100, in=180]
     (-2.5,-6) to[out=0, in=-110, looseness=0.8]
     (0.1,-4.7) to[out=70, in=-150]
     (2.9,-2.9);
    \begin{scope}[scale=0.9, rotate=50, xshift=-6.5cm, yshift=0.5cm]
     \genus
    \end{scope}
    \begin{scope}[scale=0.9, rotate=50, xshift=-6.5cm, yshift=-1.5cm]
     \genus
    \end{scope}
   \end{scope}
   
   \begin{scope}[xshift=18cm]
    \draw[very thick, color=gray!50] (-2.95,-3.3) to[out=85,in=180]
    (-2.2,-2.7) to[out=0,in=100]
    (-1.1,-3.7) to[out=-70,in=90,looseness=0.5] (-1,-5);
    \hole
    \draw[thick] (-2,1.6) -- ++(-1,-5);
    \draw[thick] (0,0) -- (-1,-5);
    \draw[very thick] (-1,-5) to[out=180,in=-95] (-3,-3.4);
    \draw[fill] (-1,-5) circle (2pt);
   \end{scope}
  \end{tikzpicture}
  \caption{The three types of surfaces with boundary: figure eight type, two holes type, cylinder (from left to right).}
  \label{fig:types_subsurfaces}
 \end{center}
\end{figure}
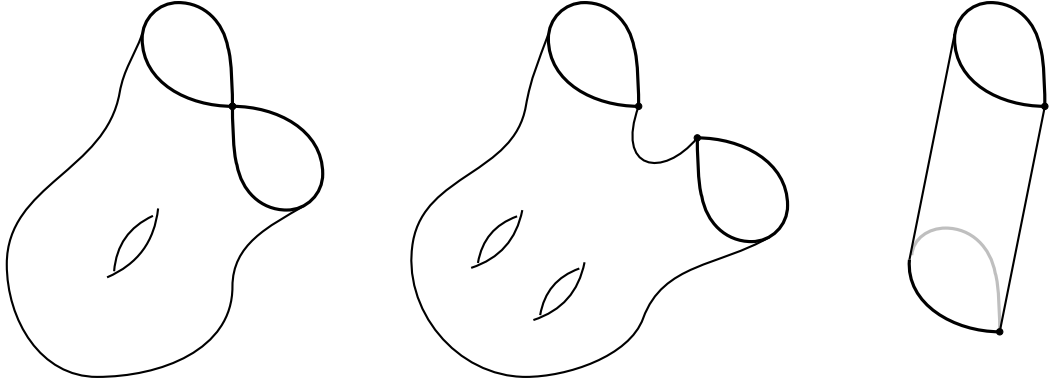

\begin{itemize}
 \item figure eight type: this has one singularity and one boundary component which consists of two saddle connections, it can have any genus $g \geq 1$
 \item two holes type: this has two boundary components that consist each of one saddle connection, it has two singularities and can have any genus $g \geq 1$
 \item cylinder: this is a special case of the two holes type where we have two boundary components but no genus
\end{itemize}

Choosing any sequence of such pieces and gluing them together in a cyclic way, gives 
us a new surface. Note that it is not possible to use only figure eight type surfaces without 
at least one surface of another type. If we would do so, the result would be a surface with two distinct points glued to each other -- which is not possible for a surface.
Also, we are in the special situation that the obtained surfaces should have only one singularity. Because of the cyclic gluing, we can have only one surface with two singularities in the sequence, 
that is, either a two hole type or a cylinder.

\subsection{Upper bound for Siegel--Veech constants for a given configuration 
with one cylinder in the minimal stratum} \label{subsec:C}

In this section, we assume that we have a cylinder, so there is no surface with two holes. We will then recover Formula 13.1 from \cite[page 142]{eskin_masur_zorich_03} but posing additionally the condition on the area. The case of a saddle connection that does not bound a cylinder and hence where there is a surface with two holes, will be discussed later in \autoref{subsec:loops}.

We first need to deal with each combinatorial type separately. 
Recall that the minimal stratum $\cH_1(2g-2)$ has three connected components which we denote by $\cH^j$ for $j=1,2,3$ from now on. The Siegel--Veech constant for a given connected component $\cH^j$ will be denoted by $c^j(V)$. We study a specific combinatorial type, described as follows. Let $p$ be the multiplicity of the saddle connection in its homology class; that is, the number of pieces in which we cut the surface (not counting the cylinder as the two boundary components of the cylinder are not both counted). Hence, we have that the surface is divided into $p$ surfaces with boundary $(X_i, \omega_i)$ of genera $g_i\geq 1$ with $\sum_{i=1}^p g_i = g-1$ and additionally a cylinder. Suppose that the first surface in the cyclic order is the cylinder.
For the other $p$ surfaces, let $a_i = 2g_i-2$ be the order of the corresponding singularity, that is, the singularity on $(X_i,\omega_i)$ has a cone angle of $(a_i+1) \cdot 2\pi$. As these $p$ surfaces are of figure eight type, we can also choose the cone angle $(a_i'+1)\cdot 2\pi$ on one side of the figure eight and have the cone angle $(a_i''+1)\cdot 2\pi$ with $a_i = a_i' + a_i''$ on the other side. Note that there are $a_i+1$ possibilities for the ordered pair $(a_i', a_i'')$.

For convenience of notation, we will refer to the datum $(p, ((a_1', a_1''),\ldots,(a_p', a_p'')))$ by $(p, \cC)$. That is, a saddle connection has the combinatorial type $(p, \cC)$ if the $p$ saddle connections in its homology class decompose the surface $(X,\omega)$ into a cylinder and $p$ surfaces of figure eight type with singularities of order $a_1,\ldots, a_p$ where the cone angles are distributed to both sides of the figure eight as given by $(a_1', a_1''), \ldots, (a_p', a_p'')$.
Let $V_{\mathrm{cyl}}(p,\cC,A)$ be the set of holonomy vectors of saddle connections of combinatorial type $(p, \cC)$ where the corresponding cylinder has area at least~$A\in [0,1)$.
Respectively, there is a Siegel--Veech constant $c^j(V_{\mathrm{cyl}}(p,\cC,A))$ on the connected component~$\cH^j$ which corresponds to counting saddle connections of combinatorial type~$(p,\cC)$ where the cylinder has area at least $A$.
We proceed to find this $c^j(V_{\mathrm{cyl}}(p,\cC,A))$ by mimicking the calculations from \cite[Section~13]{eskin_masur_zorich_03}.

For this, let $\gamma$ be the holonomy vector of a saddle connection of combinatorial type $(p, \cC)$ and~$h$ be the height of the cylinder that $\gamma$ bounds.
Note that the surface $(X,\omega)$ that we build can have any area less than or equal to $1$.
However, it can be rescaled to have area $1$ with a cylinder whose circumference is bounded from above by $\delta$ and the area is bounded from below by $A$. Hence, we have
\begin{equation*}
 |\gamma| \leq \delta \cdot \sqrt{\area (X,\omega))}
 \qquad \text{and} \qquad
 h \cdot |\gamma|\geq A \cdot \area(X,\omega).
\end{equation*}
As $\area (X,\omega) = \sum_{i=1}^p \area (X_i,\omega_i)+h|\gamma|$, the first inequality is equivalent to
\begin{equation}
 \label{eq:nontrivial}
 h \geq \frac{|\gamma|}{\delta^2}-\frac{\sum_{i=1}^p \area(X_i,\omega_i)}{|\gamma|}
\end{equation}
whereas the second inequality is equivalent to
\begin{equation}
 \label{eq:condition_on_h_from_cylinder_area}
 h\geq \frac{A}{1-A} \cdot \frac{\sum_{i=1}^p \area(X_i,\omega_i)}{|\gamma|}
 .
\end{equation}

Both of the lower bounds on $h$ have to be fulfilled and together they are sufficient to obtain a translation surface by the construction given in \cite[Section~13]{eskin_masur_zorich_03}. We distinguish two cases now, depending on whether the bound from \eqref{eq:nontrivial} or the bound from \eqref{eq:condition_on_h_from_cylinder_area} is larger and implies the other inequality.
Note~that
\begin{equation*}
 \frac{|\gamma|}{\delta^2}-\frac{\sum_{i=1}^p \area (X_i,\omega_i)}{|\gamma|}
 \geq \frac{A}{1-A} \cdot \frac{\sum_{i=1}^p \area (X_i,\omega_i)}{|\gamma|}
 ,
\end{equation*}
is equivalent to 
\begin{equation*}
 |\gamma|\geq \delta \sqrt{ \frac{\sum_{i=1}^p \area (X_i,\omega_i)}{1-A} }
 .
\end{equation*}

For the following calculation, let $d=4g$ be the real dimension of $\cH(2g-2)$ and $d_i = 4g_i$ be the real dimension of $\cH(a_i)$. Note that we have $\sum_{i=1}^p d_i = d-4$.
Furthermore, we set 
\[
r_i = \sqrt{\area(X_i, \omega_i)}
\qquad\text{and}\qquad
D(z) = \left\{ (r_1,\ldots,r_p) : \sum_{i=1}^p r_i^2 \leq z \right\}.
\]
Define $\Htc(\delta,\cC,A)$ to be the set of translation surfaces in $\Htc(\delta, A)$ 
where the saddle connection of length at most $\delta$ that bounds a cylinder of area 
$A$ has the combinatorial type $(p,\cC)$. Similar to 
\cite[Section 13.1]{eskin_masur_zorich_03}, we have 
\begin{align*}
 \nu&(\Htc(\delta,\cC,A)) = \\
 & = WM  \cdot 
  \Bigg(
 \int_{D(1-A)} \prod_{i=1}^p r_i^{d_i-1} \,dr_i 
  \Bigg(
  \int_{|\gamma| \leq \delta \sqrt{\frac {\sum r_i^2}{1-A}}} \int_{\frac{A}{1-A}\frac{\sum r_i^2}{|\gamma|} \leq h \leq \frac{1-\sum r_i^2}{|\gamma|}} \int_{0}^{|\gamma|} \,dt\,dh\,d\gamma +
  \\
  & \qquad \qquad \qquad \qquad \qquad +
  \int_{\delta \sqrt{\frac {\sum r_i^2}{1-A}} \leq |\gamma| \leq \delta} \int_{\frac{|\gamma|}{\delta^2}-\frac{\sum r_i^2}{|\gamma|} \leq h \leq \frac{1-\sum r_i^2}{|\gamma|}}\int_0^{|\gamma|}\,dt\,dh\,d\gamma
  \Bigg) \Bigg)  + o(\delta^2)
\end{align*}
where $t$ is a twist parameter of the cylinder, $W = \prod \nu(\cH(a_i))$,  $M$ is the combinatorial constant that counts how many different surfaces of area $1$ can be built for the fixed data $\delta$, $(p, \cC)$, and $A$, and $o(\delta^2)$ refers to a term that goes to zero faster than $\delta^2$ as $\delta$ goes to zero.

We first integrate over $t$ and $h$ and then over $\gamma$ to obtain
\begin{align*}
  \nu(\Htc&(\delta,\cC,A)) \\
 = \ & W M \cdot \int_{D(1-A)}\prod_{i=1}^p r_i^{d_i-1}\Bigg(
 \int_{|\gamma| \leq \delta \sqrt{\frac {\sum r_i^2}{1-A}}} \left( 1-\frac{1}{1-A} \sum r_i^2 \right) \, d\gamma \\
 & \qquad\qquad\qquad\qquad\qquad + \int_{\delta \sqrt{\frac {\sum r_i^2}{1-A}} \leq |\gamma| \leq \delta} \left(1-\frac{|\gamma|^2}{\delta^2} \right) \, d\gamma \Bigg) \prod_{i=1}^p dr_i + o(\delta^2) \\
 = \ & W M\cdot \delta^2\cdot \pi\cdot \int_{D(1-A)}\prod_{i=1}^p r_i^{d_i-1}
 \left(1-\frac{1}{1-A}\sum r_i^2\right)\left(\frac{1}{1-A} \sum r_i^2\right)\prod_{i=1}^p dr_i \\
 & \qquad \qquad + W M\cdot \delta^2\cdot \pi\cdot \int_{D(1-A)}\prod_{i=1}^p r_i^{d_i-1} 
 \left(1-\frac{1}{1-A}\sum r_i^2\right) \prod_{i=1}^p dr_i \\
 & \qquad \qquad - W M \cdot \pi\cdot \int_{D(1-A)} \prod_{i=1}^p r_i^{d_i-1}
 \left( \int_{\delta\sqrt{\frac{\sum r_i^2}{1-A}}}^\delta \, \frac{s^2}{\delta^2} \cdot 2s\,ds \right) \prod_{i=1}^p dr_i +o(\delta^2)
 .
\end{align*}

With the substitution $u=\frac{s^2}{\delta^2}$, we can write the inner most integral in the last summand~as
$$\delta^2 \cdot \int_{\frac{1}{1-A}\sum r_i^2}^{1}u \,du ,$$
giving 
\begin{align*}
  \nu(&\Htc(\delta,\cC,A))= \\
 = & W M\cdot \delta^2\cdot \pi \cdot \int_{D(1-A)}\prod_{i=1}^p r_i^{d_i-1}
 \left(1-\frac{1}{1-A}\sum r_i^2\right)\left(1 + \frac{1}{1-A} \sum r_i^2\right)\prod_{i=1}^p dr_i \\
 & - W M \cdot \delta^2 \cdot \pi \cdot \int_{D(1-A)} \prod_{i=1}^p r_i^{d_i-1}
 \left(\frac{1}{2} - \frac{1}{2(1-A)^2} \left( \sum r_i^2 \right)^2 \right) \prod_{i=1}^p dr_i +o(\delta^2) \\
 = & \pi \delta^2 \cdot W M \cdot \int_{D(1-A)}\prod_{i=1}^p r_i^{d_i-1} \\
 & \cdot \left( 1 - \frac{1}{(1-A)^2} \left( \sum r_i^2 \right)^2 - \frac{1}{2} + \frac{1}{2(1-A)^2} \left( \sum r_i^2 \right)^2 \right) \prod_{i=1}^p dr_i +o(\delta^2) \\
 = & \pi \delta^2\cdot W M \cdot \int_{D(1-A)}\prod_{i=1}^p r_i^{d_i-1}
 \left( \frac{1}{2}-\frac{1}{2}\left(\sum \left(\frac{r_i}{\sqrt{1-A}} \right)^2\right)^2 \right) \prod dr_i +o(\delta^2)\\
 = & \pi \delta^2\cdot W M \cdot \int_{D(1-A)}\prod_{i=1}^p r_i^{d_i-1} \\
 & \cdot \left( \left( 1-\sum\left(\frac{r_i}{\sqrt{1-A}}\right)^{2}\right)-\frac{1}{2}\left(1-\sum \left(\frac{r_i}{\sqrt{1-A}}\right)^{2}\right)^{2} \right) \prod dr_i+o(\delta^2) .
\end{align*}
Now make the substitution $s_i=\frac{r_i}{\sqrt{1-a}}$ and since $\sum_{i=1}^p (d_i-1) =d-4-p$ and $d=4g$, we~get
\begin{align*}
  \nu(\Htc&(\delta,\cC,A)) \\
 = &(1-A)^{\frac{p}{2}} \cdot \prod_{i=1}^p(1-A)^{\frac{(d_i-1)}{2}}\cdot \pi\delta^2 \cdot W M \\
 & \qquad \qquad \cdot \int_{D(1)}\prod_{i=1}^p s_i^{d_i-1} \left(\left(1-\sum s_i^2\right)-\frac{1}{2}\left(1-\sum s_i^2\right)^2\right) \prod ds_i \\
 = & (1-A)^{2g-2} \cdot \pi\delta^2 \cdot W M \\
 & \qquad \qquad \cdot \int_{D(1)}\prod_{i=1}^p s_i^{d_i-1}
 \left(\left(1-\sum s_i^2\right)-\frac{1}{2}\left(1-\sum s_i^2\right)^2 \right) \prod ds_i+o(\delta^2)
 .
\end{align*}

Although its appearance is different on the first view, this formula has precisely the same content as Formula 13.1 in \cite{eskin_masur_zorich_03} for $q=1$, multiplied by the factor $(1-A)^{2g-2}$ which accounts for the area bound. This factor is the same for all combinatorial types $(p, \cC)$, in particular it is independent of $p$.
The role of the factor $(1-A)^{2g-2}$ was also shown by Vorobets in \cite[Theorem 1.8]{vorobets_2005}.

Now we find the combinatorial constant $M$ and investigate its dependence on $g$.
Here we may exactly use the formulae from \cite[Section 13.3]{eskin_masur_zorich_03}.
As we are in the minimal stratum, most of the factors in the formula are equal to $1$. In this particular case, we can use the following estimate.
\begin{equation*}
 M \leq \prod_{i=1}^p (a_i+1)
\end{equation*}

The above expression for $\nu(\Htc(\delta,\cC,A))$ can now be used to calculate
$\int_{\cH_1(2g-2)} \hat f \, d\nu$ where~$\hat f$ is the function associated to 
$V_{\mathrm{cyl}}(p,\cC,A)$ which Eskin--Masur--Zorich use to the
find the associated Siegel--Veech constant. In our setting, the above integral is
still an upper bound for $\int_{\cH^j} \hat f \, d\nu$ and we get an upper bound
for $c^j(V_{\mathrm{cyl}}(p,\cC,A))$. Following the reasoning in  
\cite[page 135]{eskin_masur_zorich_03}, we have
\begin{align*}
 c^j(V_{\mathrm{cyl}}(p,\cC,A))
 \leq & M \cdot (1-A)^{2g-2} \cdot \frac{1}{2^{p-1}} \cdot \frac{\prod_i (\frac{d_i}{2}-1)!}{(\frac{d}{2}-2)!} \cdot \frac {\prod_{i=1}^p \nu(\cH(a_i))}{\nu(\cH^j)} \\
 \leq & (1-A)^{2g-2} \cdot \frac{1}{2^{p-1}} \cdot \prod_{i=1}^p (a_i+1) \cdot  \frac{\prod_i (\frac{d_i}{2}-1)!}{(\frac{d}{2}-2)!} \cdot \frac {\prod_{i=1}^p \nu(\cH(a_i))}{\nu(\cH^j)}
 .
\end{align*}
In this formula, the term $\prod_{i=1}^p \nu(\cH(a_i))$ is the previous $W$ and the factor $\frac{(1-A)^{2g-2}}{2^{p-1}} \cdot \frac{\prod_i (\frac{d_i}{2}-1)!}{(\frac{d}{2}-2)!}$ comes from the integral.
We now replace the exact values $a_i=2g_i-2$, $d_i=4g_i$, and $d=4g$.
Furthermore, we have from \cite[Theorem 1.9]{sauvaget_18}  that $\nu(\cH(a_i)) = \frac{4}{2g_i-1} \cdot (1+\cO(\frac{1}{g}))$ and $\nu(\cH_1(2g-2)) = \frac{4}{2g-1} \cdot (1+\cO(\frac{1}{g}))$.
Explicit bounds $\nu(\cH(a_i)) \leq \frac{4}{2g_i-1} \cdot (1 + \frac{2^{2^{200}}}{g_i})$ and $\nu(\cH_1(2g-2)) \leq \frac{4}{2g-1} \cdot (1 + \frac{2^{2^{200}}}{g})$ are given in the more general \cite[Theorem 1.4]{aggarwal_18}.
To avoid keeping track of the error term, we define $k\coloneqq 4(1+2^{2^{200}})$ and use the bounds $\nu(\cH(a_i)) \leq \frac{k}{2g_i-1}$ and $\nu(\cH_1(2g-2)) \leq \frac{k}{2g-1}$. Then we have
\begin{align*}
 c^j(V_{\mathrm{cyl}}(p,\cC,A))
 \leq & (1-A)^{2g-2} \cdot \frac{1}{2^{p-1}} \cdot \prod_{i=1}^p (2g_i-1) \cdot \frac{\prod_{i=1}^p (2g_i-1)!}{(2g-2)!} \cdot \frac {k^p}{\nu(\cH^j)\cdot \prod_{i=1}^p (2g_i-1)} \\
 \leq & k^p \cdot (1-A)^{2g-2} \cdot \frac{\prod_{i=1}^p (2g_i-1)!}{(2g-2)!} \cdot \frac{1}{\nu(\cH^j)}
 .
\end{align*}

\subsection{Counting of cylinders in the minimal stratum}
\label{subsec:bounds_Siegel_Veech_cylinders}

Now we want to calculate upper bounds for the Siegel--Veech 
constant for saddle connections that bound a cylinder. For a given multiplicity $p$, we 
consider all ways how to decompose the surface into a cylinder and $p$ surfaces with boundary and 
how to decompose the angle of the singularity of order $a_i$ in each of the surfaces.
That is, for a connected component $\cH^j$, we define $V_{\mathrm{cyl}}(p, A)$ to be the set of 
holonomy vectors for all saddle connections bounding a cylinder of area at least $A$
with multiplicity $p$. Let $c^j(V_{\mathrm{cyl}}(p,A))$ be the associated 
Siegel--Veech constant. Then, $c^j(V_{\mathrm{cyl}}(p,A))$ is calculated as 
sum over all possible combinatorial types $(p,\cC) = (p, ((a_1',a_1''),\ldots,(a_p',a_p'')))$ for fixed~$p$:
\begin{align*}
 c^j(V_{\mathrm{cyl}}(p,A)) & = \sum_{\cC} c^j\big(V_{\mathrm{cyl}}(p,\cC, A)\big)\\
& \leq \sum_{g_1+\ldots+g_p = g-1} \prod_{i=1}^p (a_i+1) \cdot k^p \cdot (1-A)^{2g-2} \cdot \frac{\prod_{i=1}^p (2g_i-1)!}{(2g-2)! \cdot \nu(\cH^j)} \\
 & \leq \frac{k^p \cdot (1-A)^{2g-2}}{\nu(\cH^j)} \cdot \sum_{g_1+\ldots+g_p = g-1} \frac{\prod_{i=1}^p (2g_i)!}{(2g-2)!} 
\end{align*}
We use $a_i + 1 \leq 2g_i$ here to obtain the third line.

Note that $\prod_{i=1}^p (2g_i)!$ and $(2g-2)!$ have the same number of factors. We have that $\prod_{i=1}^p (2g_i)!$ is the largest when all but one $g_i$ are equal to $1$. In this situation, let $g_1$ be the largest one, that is $g_1 = g-1 - (p-1)$. Therefore, we have
\begin{equation*}
 \prod_{i=1}^p (2g_i)! \leq 2^{p-1} \cdot (2g-2p)!
\end{equation*}
which implies
\begin{align*}
 \frac{\prod_{i=1}^p (2g_i)!}{(2g-2)!}
 \leq & 2^{p-1} \cdot \frac{(2g-2p)!}{(2g-2)!} \\
 = & \frac{2^{p-1}}{(2g-2)(2g-3)\cdot \cdots \cdot (2g-2p+1)} \\
 \leq & \frac{1}{(2g-2)(2g-3)\cdot \cdots \cdot (2g-p)}
\end{align*}
for $p\geq 2$ and $\frac{(2g_1)!}{(2g-2)!} = 1$ for $p = 1$.

This upper bound is independent of the choice of the $g_i$. There are $\binom{g-2}{p-1}$ choices of the~$g_i$.
This is because we can consider an ordered set with $g-1$ elements and divide it into $p$ subsets of cardinality $\geq 1$ by specifying which elements are the last in their corresponding subsets. The last one in the whole ordered set has to be the last one of a subset. Apart from that, we can choose any $p-1$ elements out of the remaining $g-2$ elements to be last ones.

Therefore, we can calculate:
\begin{align*}
 \sum_{g_1+\ldots+g_p = g-1} \frac{\prod_{i=1}^p (2g_i)!}{(2g-2)!}
 & \leq \binom{g-2}{p-1} \cdot \frac{1}{(2g-2)(2g-3) \cdots  (2g-p)} \\
 & = \frac{1}{(p-1)!} \cdot \frac{(g-2)(g-3) \cdots (g-p)}{(2g-2)(2g-3) \cdots (2g-p)} \\
 & \leq \frac{1}{(p-1)!}
\end{align*}
for $p\geq 2$ and the same upper bound is true for $p=1$.

Still fixing a connected component $\cH^j$, we define $V_{\mathrm{cyl}}(A)$ to be the set of 
holonomy vectors for all saddle connections bounding a cylinder of area at least $A$. Let $c^j(V_{\mathrm{cyl}}(A))$ be the associated Siegel--Veech constant. 
To compute $c^j(V_{\mathrm{cyl}}(A))$, we have to sum over all $p$. This gives us:
\begin{align}
 c^j(V_{\mathrm{cyl}}(A))
 = & \sum_{p=1}^{g-1} c^j\big(V_{\mathrm{cyl}}(p,A)\big) \notag \\
 \leq & \frac{(1-A)^{2g-2}}{\nu(\cH^j)} \cdot
    \sum_{p=1}^{g-1} k^p \cdot \frac{1}{(p-1)!} \notag \\
 \leq & \frac{(1-A)^{2g-2}}{\nu(\cH^j)} \cdot k \cdot \sum_{p=0}^{g-2} \frac{k^p}{p!} \notag \\
 \leq & \frac{k \cdot e^{k} \cdot (1-A)^{2g-2}}{\nu(\cH^j)}
      .  \label{eq:cyl-estimate}
\end{align}

We are now ready to prove \autoref{thm:measure_cylinders_large_area} for the minimal stratum.
\begin{proof}[Proof of \autoref{thm:measure_cylinders_large_area}]
Let $V_{\mathrm{cyl}}(A)$ be defined as above, $f \colon \mathbb{R}^2\to \mathbb R$ be the characteristic function of the ball of radius $\delta$ and 
$\hat f \colon\, \cH^j \to \mathbb{R}$ be the associated function defined in 
\autoref{Eq:f-hat}. 

We argue that $\hat{f}(X) = N_{\mathrm{cyl}}(X, \delta, A)$ outside of 
a measure zero set.
The two numbers are not equal if and only if there exists more than one cylinder that is bounded by a given saddle connection in $V_{\mathrm{cyl}}(A)$ or if there is a cylinder of the desired type with more than one saddle connection on a boundary component.
Recall that the measure on the stratum is given as Lebesgue measure in period coordinates and
in the minimal stratum, the relative homology is the same as the absolute homology.
Consider the set $\cH_{\mathrm{parallel}}$ of translation surfaces that have two non-homologous 
saddle connections whose holonomy vectors have the same direction. As $\cH_{\mathrm{parallel}}$ can be defined locally in period coordinates, it is a lower-dimensional subset of $\cH^j$ and hence has measure zero. This has two consequences. First, there is only a measure zero set of translation surfaces that has a cylinder with a boundary component that contains more than
one saddle connections. That is, generically, for every cylinder of area $A$, 
we have a vector in $V_{\mathrm{cyl}}(A)$. Secondly, the set of translation surfaces 
that have more than one cylinder giving the same vector in $V_{\mathrm{cyl}}(A)$
has also measure~zero. 

We also have
 \begin{equation*}
   \int_{\cH^j} \hat{f} \, d\nu 
   = c^j(V_{\mathrm{cyl}}(A)) \cdot \int_{\mathbb {R}^2} f \, dxdy \cdot  \nu(\cH^j)
   \leq k e^k \cdot \pi \,\delta^2 \, (1-A)^{2g-2}.
 \end{equation*}
We can now estimate the desired integral by using the bounds on the integrals for all three connected components $\cH^j$ on the stratum.
 \begin{align*}
  \frac{\int_{\cH_1(2g-2)} N_{\mathrm{cyl}}(X, \delta, A) \, d\nu (X)}{\nu\big( \cH_1(2g-2) \big)}
  \leq & \frac{1}{\nu\big( \cH_1(2g-2) \big)} \cdot 
      \sum_{i=1}^3 \int_{\cH^j} \hat{f} \, d\nu  \\
  \leq &  \frac{2g}{3.9}  \cdot \sum_{i=1}^3 ke^k \cdot \pi \,\delta^2 \, (1-A)^{2g-2}\\
  \leq & 2 k e^k \cdot \, g \, \delta^2 \, (1-A)^{2g-2}
  .
 \end{align*}
 
 For the second to last line, we use again \cite[Theorem 1.4]{aggarwal_18}, this time in the form that $\nu\big( \cH_1(2g-2) \big) \geq \frac{3.9}{2g}$ for large values of $g$.
 
 The estimate for $\nu(\Htc(\delta, A))$ follows from the fact that $X \in \Htc(\delta, A)$ if and only if $N_{\mathrm{cyl}}(X, \delta, A) \geq 1$. 
\end{proof}

\subsection{Counting of saddle connections in the minimal stratum} \label{subsec:loops} 

For the sake of completeness, we now turn to the case of saddle connections that do not bound a cylinder. Let $p$ again be the multiplicity of the saddle connection. Then the surface decomposes into $p$ surfaces with boundary of which $p-1$ are of figure eight type and exactly one is of two holes type. In particular, there is no cylinder.

Suppose that the first surface is the surface of two holes type.
Let $b_1',b_1''\geq 0$ be integers such that the interior angle at the holes is $(2b_1'+3)\pi$ and $(2b_1''+3)\pi$ with $b_1'+b_1''= 2g_1 -2$.
Then the real dimension of the stratum $\cH(b_1'', b_1'')$ of the first surface is $d_1 = 4g_1+2$ and the volume of the stratum $\cH(b_1'', b_1'')$ is approximately $\frac{4}{(b_1'+1)(b_1''+1)}$ (see \cite[Theorem 1.4]{aggarwal_18}). We use again the bound $\nu(\cH(b_1'', b_1'')) \leq \frac{k}{(b_1'+1)(b_1''+1)}$ with $k = 4(1+2^{2^{200}})$.

Similarly to before, this data defines the combinatorial type $(p,\cC)$. We consider the corresponding set $V_{\mathrm{loop}}(p,\cC)$ of holonomy vectors of saddle connections of combinatorial type~$(p,\cC)$ that do not bound a cylinder.
Then Formula 13.1 from \cite{eskin_masur_zorich_03} gives us that in the situation of having exactly one surface of two holes type and no cylinder, the Siegel--Veech constant $c^j(V_{\mathrm{loop}}(p,\cC))$ for this data is bounded in the following way.
\begin{align*}
 c^j(V_{\mathrm{loop}} & (p,\cC)) \\
 & \leq \frac{1}{2^{p-1}} \cdot (b_1'+1) (b_1''+1)\cdot \prod_{i=2}^p (a_i+1) \cdot  \frac{\prod_i (\frac{d_i}{2}-1)!}{(\frac{d}{2}-2)!} \cdot \frac {\nu(\cH(b_1', b_1'')) \cdot \prod_{i=2}^p \nu(\cH(a_i))}{\nu(\cH^j)} \\
 &\leq  \frac{1}{2^{p-1}} \cdot (b_1 '+1) (b_1''+1) \\
 & \quad \cdot \prod_{i=2}^p (2g_i-1) \cdot \frac{(2g_1)! \cdot \prod_{i=2}^p (2g_i-1)!}{(2g-2)!} \cdot \frac {k^p}{(b_1'+1) (b_1''+1) \cdot \prod_{i=2}^p (2g_i-1) \cdot \nu(\cH^j)}
 \\
 & \leq k^p \cdot \frac{2g_1 \cdot \prod_{i=1}^p (2g_i-1)!}{(2g-2)!} \cdot \frac {1}{\nu(\cH^j)}
 .
\end{align*}

Note that the bound for this Siegel--Veech constant $c^j(V_{\mathrm{loop}}(p,\cC))$ differs exactly by a factor of $2g_1$ from the bound for the Siegel--Veech constant $c^j(V_{\mathrm{cyl}}(p,\cC, 0))$ where we have a cylinder without a condition on the area.
Hence, we can skip the calculations for a given multiplicity $p$ and deduce directly
\begin{equation} \label{eq:loop-estimate}
 c^j(V_{\mathrm{loop}}) \leq 2g_i \cdot \frac{k e^k}{\nu(\cH^j)} \leq g \cdot \frac{2k e^k}{\nu(\cH^j)}
 .
\end{equation}

Let $c^j(V_{\mathrm{sc}})$ be the Siegel--Veech constant for all saddle connections. 
Then 
\[
c^j(V_{\mathrm{sc}}) \leq c^j(V_{\mathrm{loop}})  + c^j(V_{\mathrm{cyl}}(0)).  
\]
Combining the estimates in \autoref{eq:cyl-estimate} for $A=0$ and 
\autoref{eq:loop-estimate} we get
\[
c^j(V_{\mathrm{sc}}) \leq  \frac{3ke^k \cdot g}{\nu(\cH^j)}.
\]
The rest of the proof of Theorem~\ref{Thm:sc} is then completely analogous to the
proof of \autoref{thm:measure_cylinders_large_area}.

\subsection{General strata}

We now turn our attention to the general case of a stratum $\cH(\kappa)$.
The proof follows the same steps as for the minimal stratum.
We will outline here how to combine the arguments from the previous sections of the appendix with the results of \cite{eskin_masur_zorich_03, vorobets_2005, aggarwal_18}. We also want to refer to the newer article \cite{chen_moeller_sauvaget_zagier_19} from which more direct deductions are possible, in particular for disconnected strata.

We start with the case of cylinders.
For this, fix a connected component $\cH$ of a stratum~$\cH(\kappa)$.
By applying the volume estimates from \cite{aggarwal_18} to the formulae from \cite{eskin_masur_zorich_03}, Zorich in an appendix to \cite{aggarwal_18} gets the following estimates for connected strata.
For cylinders where one boundary component is a saddle connection through a singularity of order $k_1$ and the other through a distinct singularity of order $k_2$ and where the multiplicity is~$1$, the Siegel--Veech constant is (up to lower order terms) 
$\frac{(k_1+1)(k_2+1)}{d-2}$ for large $g$.
If the singularities are the same on both boundary components, we have that the Siegel--Veech constant is (again, up to lower order terms)  $\frac{1}{2} \cdot \frac{(k_1+1)(k_1-1)}{d-2}$ for large $g$.
Note that in the case of non-connected strata, these are not estimates for the Siegel--Veech constants on the stratum but we have the upper bounds
$\frac{(k_1+1)(k_2+1)}{d-2} \cdot \frac{\nu(\cH_1(\kappa))}{\nu(\cH_1)}$ and  $\frac{1}{2} \cdot \frac{(k_1+1)(k_1-1)}{d-2} \cdot \frac{\nu(\cH_1(\kappa))}{\nu(\cH_1)}$
 on the Siegel--Veech constants on all of the connected components.

We now have to compare the Siegel--Veech constants for multiplicity~$1$ with these for higher multiplicity.
In \cite{aggarwal_18_SV}, Aggarwal shows that the Siegel--Veech constants for saddle connections of higher multiplicity are of lower order than the ones for saddle connections of multiplicity~$1$. As the terms for the saddle connection Siegel--Veech constants differ from the terms for the cylinder Siegel--Veech constants by a factor of $d-2$, the same combinatorial arguments hold for the cylinder Siegel--Veech constants.
Hence, the Siegel--Veech constants for cylinders with restrictions on the order of the singularity but without the restriction on the multiplicity are bounded by $K \cdot \frac{(k_1+1)(k_2+1)}{d-2}\cdot \frac{\nu(\cH_1(\kappa))}{\nu(\cH_1)}$ and $K \cdot \frac{(k_1+1)(k_1-1)}{d-2} \cdot \frac{\nu(\cH_1(\kappa))}{\nu(\cH_1)}$, respectively, for large values of $g$ and some constant $K>0$.

Recall that $V_\mathrm{cyl} \coloneqq V_\mathrm{cyl}(0)$ is the set of holonomy vectors for all saddle connections bounding a cylinder.
To obtain a bound on the Siegel--Veech constant $c(V_\mathrm{cyl})$ for $\cH$, we have to sum over all possible ordered pairs $(k_1, k_2)$:
\begin{align*}
 c(V_\mathrm{cyl})
 & \leq K \cdot \frac{1}{d-2} \cdot \frac{\nu(\cH_1(\kappa))}{\nu(\cH_1)} \cdot
 \Bigg( \sum_{i\neq j}(k_i+1)(k_j+1) + \sum_{k_i\geq 2} (k_i+1)(k_i-1)\Bigg) \\
 & \leq K \cdot \frac{1}{d-2} \cdot \frac{\nu(\cH_1(\kappa))}{\nu(\cH_1)} \cdot \sum_{i,j} (k_i+1)(k_j+1) \\
 & = K \cdot \frac{1}{d-2} \cdot \frac{\nu(\cH_1(\kappa))}{\nu(\cH_1)} \cdot \left(\sum_i (k_i+1)\right)^2 \\
 & \leq K \cdot \frac{1}{d-2}(4g-4)^2 \cdot \frac{\nu(\cH_1(\kappa))}{\nu(\cH_1)}
\end{align*}
In the last inequality, we used $\sum_i (k_i + 1) \leq \sum_i k_i + \ell \leq (2g-2) + (2g-2)$. Note
that $d-2$ is larger than $2g-2$ and hence we have
\begin{align*}
 c(V_\mathrm{cyl})
  &\leq K \cdot \frac{1}{d-2}(4g-4)^2 \cdot \frac{\nu(\cH_1(\kappa))}{\nu(\cH_1)}\\
  &\leq K \cdot \frac{1}{2g-2}(4g-4)^2 \cdot \frac{\nu(\cH_1(\kappa))}{\nu(\cH_1)}\\
  &\leq 2K \cdot (4g-4)\cdot \frac{\nu(\cH_1(\kappa))}{\nu(\cH_1)}
   \leq 8K \cdot g \cdot \frac{\nu(\cH_1(\kappa))}{\nu(\cH_1)} .
\end{align*}

As in the proof of \autoref{thm:measure_cylinders_large_area}, we let $f$ be the characteristic function of the ball of radius~$\delta$ and~$\hat f \colon\, \cH^j \to \mathbb{R}$ be the associated function defined in \autoref{Eq:f-hat}. 
Again, $\hat{f}(X) = N_{\mathrm{cyl}}(X, \delta)$ outside of a measure zero set.
This is because, otherwise the holonomy vectors associated to two not homologous 
saddle connections are parallel and this is a measure zero property. 
As every stratum has at most three connected components, we can again use the calculation
 \begin{align*}
  \frac{\int_{\cH_1(\kappa)} N_{\mathrm{cyl}}(X, \delta) \, d\nu (X)}{\nu\big( \cH_1(\kappa) \big)}
  \leq & \, \frac{1}{\nu\big( \cH_1(\kappa) \big)} \cdot 3 \cdot \int_{\cH} \hat{f} \, d\nu  \\
  \leq & \, \frac{1}{\nu\big( \cH_1(\kappa) \big)} \cdot 3 \cdot 8K \cdot g \cdot \frac{\nu(\cH_1(\kappa))}{\nu(\cH_1)} \cdot \pi \delta^2 \cdot \nu\big(\cH_1 \big) \\  
  \leq & \, 24 K \cdot g \, \delta^2
  .
 \end{align*}

Including the requirement on the area of the cylinder gives a factor of $(1-A)^{d-2}$ in the very first calculation of Siegel--Veech constants for a given configuration (compare the proof for the minimal stratum or \cite{vorobets_2005}).
Hence, this factor carries through the full proof and appears in the end as claimed.

For the case of saddle connections, a comparison of Corollary 1 and 3 with Corollary 4 and~5 in the appendix of \cite{aggarwal_18} shows that the estimates of the corresponding Siegel--Veech constants are larger by a factor of $d-2$ than in the case of cylinders. The inequality $3g \leq d-2$ then implies \autoref{Thm:sc}.

\bibliographystyle{amsalpha}
\bibliography{literature}

\end{document}